\documentclass[journal]{IEEEtran}
\usepackage{amsmath}
\usepackage{graphicx}
\usepackage{dsfont}
\usepackage{epstopdf}
\usepackage{subfigure}
\usepackage{algorithm}
\usepackage{algpseudocode}
\usepackage{cite}
\usepackage{multicol}
\usepackage{multirow}
\usepackage{mathtools}
\usepackage{amssymb}
\usepackage{color}
\usepackage{bm}
\usepackage{indentfirst}
\usepackage{url}
\usepackage{nomencl}
\usepackage{booktabs}
\usepackage{cite}
\usepackage[numbers,sort&compress]{natbib}
\usepackage{graphicx}
\usepackage{subfigure}
\usepackage{algorithm}
\usepackage{algpseudocode}
\usepackage{cite}
\usepackage{multicol}
\usepackage{multirow}
\usepackage{mathtools}
\usepackage{amssymb}
\usepackage{color}
\usepackage{bm}
\usepackage{indentfirst}
\usepackage{url}
\usepackage{nomencl}
\usepackage{booktabs}
\usepackage{cite}
\usepackage[numbers,sort&compress]{natbib}
\usepackage{amsthm}
\newtheorem{Theorem}{Theorem}
\newtheorem{Lemma}{Lemma}

\newtheorem{Example}{Example}
\newtheorem{Remark}{Remark}

\newtheorem{Assumption}{Assumption}
\allowdisplaybreaks

\begin{document}

\title{Distributed Optimization with Coupling Constraints in Multi-Cluster Networks Based on Dual Proximal Gradient Method} 


\author{Jianzheng Wang, Guoqiang~Hu,~\IEEEmembership{Senior Member,~IEEE} 
}
\maketitle

\begin{abstract}                          
In this work, we consider solving a distributed optimization problem in a multi-agent network with multiple clusters. In each cluster, the involved agents cooperatively optimize a separable composite function with a common decision variable. Meanwhile, a global cost function of the whole network is considered associated with an affine coupling constraint across the clusters. To solve this problem, we propose a cluster-based dual proximal gradient algorithm by resorting to the dual problem, where the global cost function is optimized when the agents in each cluster achieve an agreement on the optimal strategy and the global coupling constraint is satisfied. In addition, the proposed algorithm allows the agents to only communicate with their immediate neighbors with heterogeneous step-sizes. The computational complexity of the proposed algorithm with simple-structured cost functions is discussed and an ergodic convergence with rate $O(\frac{1}{T})$ is guaranteed ($T$ is the index of iterations).
\end{abstract}


\section{Introduction}

\subsection{Background and Motivation}

Distributed optimization problems (DOPs) have drawn much attention in the recent few years due to their wide applications in practical problems, such as task assignment in multi-robot networks \cite{luo2014provably}, machine learning problems \cite{lee2017speeding}, and economic dispatch in power systems \cite{bai2017distributed}. In those problems, each agent usually maintains an individual decision variable, and the optimal solution of the whole network is achieved through multiple rounds of local communication and decision making.

To adapt to the arising computational burden and privacy preservation issues in the practical implementations with large-scale data sets, DOPs with multiple agent clusters were discussed recently \cite{guo2017distributed}. Generally speaking, a multi-cluster network corresponds to the conventional multi-agent network where each node is broken into a cluster of sub-nodes such that the overall computation task for the node can be separated effectively \cite{shi2019multi}. In a multi-cluster-based DOP, to achieve the respective objective of each cluster and the whole network, the state updates and information exchanges are usually analyzed at two levels: cluster level and system level, which brings more complications into the algorithm development than the single-level counterparts. In addition, existing research works in this field expose a significant gap in considering coupling constraints across the clusters, which limits their usages in various coupling constrained problems \cite{zhou2017multi,notarnicola2016duality,chang2016proximal}. Therefore, in this work, we aim to develop a distributed optimization algorithm for multi-cluster networks with general affine coupling constraints.

\subsection{Literature Review}

Existing research works on DOPs usually assume that the communication graph of multi-agent networks is incomplete due to the inadequate communication infrastructures or communication failures among the agents. To tackle this issue, a valid alternative is applying graph theory in modeling the communication links, by which each agent only needs to access the state information from its immediate neighbors \cite{pang2019randomized,ning2018distributed,wang2019distributed}. In addition, to deal with the coupling constraints, consensus protocol was applied in \cite{chang2014distributed,simonetto2016primal}, where certain agreement on the state of the network can be achieved only through local communications. To further minimize the requirement on the acquisition of the global information, some more distributed algorithms were proposed by \cite{mosk2010fully,necoara2017fully,wang2021distributed}, where the step-size of the agents can be determined by local information. Some DOPs dealing with time-varying networks were discussed in \cite{falsone2017dual,yang2010distributed}. Some distributed algorithms with accelerated convergence rates can be referred to in \cite{necoara2015linear,alghunaim2019proximal}.

However, when we extend the conventional multi-agent networks to the multi-cluster cases, the algorithms in the aforementioned works fail to illuminate the realization mechanisms within and across the clusters. Regarding this, an unconstrained DOP with multiple clusters was investigated in \cite{guo2017distributed}, where the agreement of the agents in each cluster is achieved by the consensus protocol. Then, the optimal solution of the whole network is achieved when the leader of a cluster executes update in a sequential order by communicating with the leader of neighbor clusters. With a similar network setting, the authors of \cite{li2019gossip} further proposed a gossip-based communication mechanism among the leaders with a supermartingale convergence guarantee. Then, a random-sleep updating strategy was proposed in \cite{shi2019multi}, where the followers in each cluster can be in an active or inactive mode based on Bernoulli distribution.

In contrast with the existing works with similar problem set-ups, the new features of this work are summarized as follows.
\begin{itemize}
  \item We newly consider a DOP in a multi-cluster network with affine coupling constraints across the clusters, which has a significant gap in both algorithm development and application potential from the cluster-based unconstrained DOPs studied in \cite{guo2017distributed,li2019gossip,shi2019multi}. In addition, the proposed multi-cluster network model is more distributed than those in \cite{guo2017distributed,li2019gossip,shi2019multi} in the sense that we do not set any leader agent who is designed to access the information of other leaders directly.
  \item Our considered cluster-based DOPs involve two levels of optimization objectives: intra-cluster consensus and extra-cluster optimization. Specifically, the agreement of each cluster is achieved by the consensus protocol among the agents involved. Then, the optimal solution of the whole network is achieved with the presence of the coupling constraint. Therefore, it can be technically challenging to directly apply the existing algorithms without clusters, as studied in \cite{chang2014distributed,li2020distributed,notarnicola2017duality,wang2021distributed,yang2010distributed,falsone2017dual,simonetto2016primal,mosk2010fully,necoara2017fully}, due to the heterogeneity of the optimization objectives within and across the clusters.
  \item To solve the DOP of interest, we propose a cluster-based dual proximal gradient (CDPG) algorithm. As an advantage of the algorithm, if the cost functions are with some simple structures (can be non-smooth), one only need to update the dual variables by some simple operations without any costly inner-loop optimizations, which leads to the reduction of the overall computational complexity. This benefit is closely based on the dual proximal gradient method, as also studied in \cite{notarnicola2016asynchronous,beck2014fast,kim2016fast,wang2021composite,wang2021distributed}. But differently, we focus on solving DOPs with general affine coupling constraints, which are not considered in \cite{notarnicola2016asynchronous,beck2014fast,kim2016fast,wang2021composite}. In addition, the proposed CDPG algorithm can avoid the bias error of a near-optimal convergence in \cite{wang2021distributed} with the same convergence rate guarantee.
\end{itemize}

The overall contributions of this work are summarized as follows.

\begin{itemize}
  \item We consider a class of DOPs in multi-cluster networks with composite cost functions and general affine coupling constraints. In this problem, the decision of each cluster is based on the consensus of the agents involved. The optimal solution of the whole network is achieved when all the agents only communicate with their immediate neighbors.
  \item A CDPG algorithm is proposed by resorting to the dual problem. Provided that the cost functions in the primal problem are simple-structured, the overall computational complexity can be reduced without any costly inner-loop optimization. In addition, the agents can update with heterogeneous step-sizes and an ergodic convergence rate $O(\frac{1}{T})$ can be guaranteed. The feasibility of the proposed algorithm is verified by solving a social welfare optimization problem in commodity market and an economic emission dispatch problem in energy market in the simulation.
\end{itemize}

The rest of this paper is organized as follows. Section \ref{sa2} presents some frequently used definitions and properties employed by this work. Section \ref{sa3} formulates the optimization problem of interest and provides the basic assumptions. Section \ref{sa4} presents the main result of this work, including the dual problem formulation, the proposed distributed optimization algorithm, and the convergence analysis. The feasibility of the proposed algorithm is numerically verified by two motivating applications in Section \ref{sa6}. Section \ref{sa7} concludes this paper.

\section{Preliminaries}\label{sa2}

Some frequently used notations, definitions, and relevant properties of proximal mapping and Fenchel conjugate are provided in this section.

\subsection{Notations}

$\mathds{N}$ and $\mathds{N}_+$ denote the non-negative and positive integer spaces, respectively. Let $\mid {A}\mid$ be the cardinality of set ${A}$. $\mathds{R}_+^n$ denotes the $n$-dimensional Euclidian space only with non-negative real elements.
Operator $(\cdot)^{\top}$ represents the transpose of a matrix.
$\prod_{n \in B} {A}_n$ and $\bigcap_{n \in B} {A}_n$ denote the Cartesian product and interaction of sets ${A}_1$ to ${A}_{|B|}$, respectively.
$\mathbf{relint}{A}$ represents the relative interior of set ${A}$. $\| \cdot\|_1$ and $\| \cdot \|$ refer to the $l_1$- and $l_2$-norms, respectively. Define $\| \mathbf{u}\|^2_{\mathbf{X}}= \mathbf{u}^{\top} \mathbf{X}\mathbf{u}$ with $\mathbf{X}$ a square matrix. $\otimes$ is Kronecker product. $\tau_{\mathrm{max}}(A)$ denotes the largest eigenvalue of a positive semi-definite matrix $A$. $\mathbf{I}_n$ is an $n$-dimensional identity matrix and $\mathbf{O}_{n \times m}$ is an $(n \times m)$-dimensional zero matrix. $\mathbf{1}_n$ and $\mathbf{0}_n$ denote the $n$-dimensional column vectors with all entries being 1 and 0, respectively. Let $\mathrm{diag}[\mathbf{X}_n]_A$ be a diagonal-like matrix with $\mathbf{X}_n$ placed on the diagonal according to index $n = 1,2,...,|A|$, i.e., $\mathrm{diag}[\mathbf{X}_n]_A = \left[
                                   \begin{array}{ccc}
                                     \mathbf{X}_1 &  & \mathbf{O} \\
                                      & \ddots &  \\
                                     \mathbf{O} &  & \mathbf{X}_{|A|} \\
                                   \end{array}
                                 \right]
$.
\subsection{Proximal Mapping}
A proximal mapping of a proper, convex, and closed function $\psi: \mathds{R}^n \rightarrow (-\infty,+\infty]$ is defined by $\mathrm{prox}^{\alpha}_{\psi} [\mathbf{v}] =  \arg \min_{\mathbf{u}} ( \psi(\mathbf{u}) + \frac{1}{2{\alpha}} \| \mathbf{u} - \mathbf{v}\|^2 )$, ${\alpha}>0$, $\mathbf{v} \in \mathds{R}^n$.
A generalized version of proximal mapping can be defined as
\begin{align}\label{13}
\mathrm{prox}^{\mathbf{X}}_{\psi} [\mathbf{v}] =  \arg \min_{\mathbf{u}} ( \psi(\mathbf{u}) + \frac{1}{2} \| \mathbf{u} - \mathbf{v}\|^2_{\mathbf{X}^{-1}} ),
\end{align}
with $\mathbf{X} \in \mathds{R}^{n \times n}$ a positive definite matrix \cite{notarnicola2016asynchronous}.

\subsection{Fenchel Conjugate}\label{de1}
$\psi: \mathds{R}^n \rightarrow (-\infty,+\infty]$ is a proper function. Then, the Fenchel conjugate of $\psi$ is defined by $\psi^{\diamond}(\mathbf{v})= \sup_{\mathbf{u}} \{\mathbf{v}^{\top} \mathbf{u}-\psi(\mathbf{u})\}$, which is convex \cite[Sec. 3.3]{borwein2010convex}.

\begin{Lemma}\label{md}
(Extended Moreau Decomposition \cite[Thm. 6.45]{beck2017first}) $\psi: \mathds{R}^n \rightarrow (-\infty,+\infty]$ is a proper, convex, and closed function. $\psi^{\diamond}$ is its Fenchel conjugate. Then for some $\mathbf{v} \in \mathds{R}^n$ and $\alpha >0$, we have
\begin{align}\label{am}
 \mathbf{v} = \alpha \mathrm{prox}^{\frac{1}{\alpha}}_{\psi^{\diamond}} [\frac{\mathbf{v}}{\alpha}] + \mathrm{prox}^{\alpha}_{\psi} [\mathbf{v}].
\end{align}
\end{Lemma}

\begin{Lemma}\label{l1}
\cite[Lemma V.7]{notarnicola2016asynchronous} Let $\psi: \mathds{R}^n \rightarrow (-\infty,+\infty]$ be a proper, closed, $\sigma$-strongly convex function and $\psi^{\diamond}$ be its Fenchel conjugate, $\sigma >0$. Then,
\begin{align}
\arg \max\limits_{\mathbf{u}} ( \mathbf{v}^{\top} \mathbf{u} - \psi(\mathbf{u})) = \nabla_{\mathbf{v}} \psi^{\diamond}(\mathbf{v})
\end{align}
and $\nabla_{\mathbf{v}} \psi^{\diamond}(\mathbf{v})$ is Lipschitz continuous with Lipschitz constant $\frac{1}{\sigma}$.
\end{Lemma}

\section{Problem Formulation}\label{sa3}

The considered network model, problem formulation, and relevant assumptions are introduced in this section.

\subsection{Network Model}
Consider a multi-agent network ${{G}}$, which is composed of cluster set $V=\{1,...,N\}$. Cluster $i \in V$ is defined by ${G}_i= \{{V}_i,{E}_i\}$ with agent set ${V}_i=\{1,...,n_i\}$ and undirected edge set $E_i =\{e_1,e_2,...,e_{|E_i|}\}  \subseteq \{(j,l)| j \in {V}_i, l\in {V}_i, j \neq l\}$ (no self-loop). Then $G$ can be described by agent set $\bar{V}=\{1,...,\sum_{l=1}^{N}n_l\}$ and undirected edge set $E = \{\bar{e}_1,\bar{e}_2,...,\bar{e}_{|E|}\} \subseteq \{(i,j)| i \in \bar{V}, j\in \bar{V}, i \neq j\}$ (no self-loop). In $\bar{V}$, the index of the $j$th agent in  cluster $i$ is relabeled by $n_{ij} = \sum_{l=0}^{i-1}n_l + j$, i.e., the agents are relabeled from cluster $1$ to $N$ according to the index in each respective cluster. Let ${{V}}_i^j =  \{ l | (j,l) \in {E}_i\}$ and $\bar{V}^k =  \{ l | (k,l) \in {E}\}$ be the neighbor sets of the $j$th agent in $G_i$ and the $k$th agent in $G$, respectively.
Let $\mathbf{L}^i \in \mathds{R}^{n_i \times n_i}$ be the Laplacian matrix of $G_i$, where the $(j,l)$th entry is defined by \cite{chung1997spectral}
\begin{align}\label{}
& [\mathbf{L}^i]_{jl}=\left\{
\begin{array}{ll}
  |V^j_i| & \hbox{if $j=l \in  V_i$} \\
  -1 & \hbox{if $(j,l) \in E_i$} \\
  0 & \hbox{otherwise}
\end{array}
\right..
\end{align}

With a given indexing protocol of the vertices in $V_i$, the index of edges can be decided as follows. For any two distinct edges $e_k = (k_1,k_2) \in E_i$ and $e_v=(v_1,v_2) \in E_i$, if $\min\{k_1,k_2\} > \min\{v_1,v_2\}$, then $k > v$, and vice versa. For the case $\min\{k_1,k_2\} = \min\{v_1,v_2\}$, if $\max\{k_1,k_2\} > \max\{v_1,v_2\}$, then $k > v$, and vice versa.

\begin{Example}\label{e1}
Based on the vertex indexing protocol in Fig. \ref{exa}, the edge indices are given by $e_1=(1,2)$, $e_2=(2,4)$, $e_3=(2,5)$, $e_4=(3,4)$, and $e_5=(4,5)$.
\begin{figure}[H]
  \centering
  \includegraphics[width=5cm]{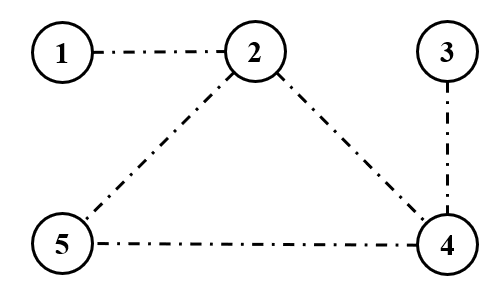}\\
  \caption{A network example.}\label{exa}
\end{figure}
\end{Example}

The graph $G_i$ can also be described by an incidence matrix $\mathbf{G}_i \in \mathds{R}^{n_i \times |E_i|}$, where the rows and columns are indexed by vertices and edges, respectively \cite{dimarogonas2010stability}. Based on the above indexing protocol of edges, the incidence matrix of $G_i$ is defined by
\begin{align}\label{}
[\mathbf{G}_i]_{jk}=\left\{
\begin{array}{ll}
  1 & \hbox{if $e_k=(j,l) \in E_i$ and $j>l$} \\
  -1 & \hbox{if $e_k=(j,l) \in E_i$ and $j<l$} \\
  0 & \hbox{otherwise}
\end{array}
\right..
\end{align}
Similarly, we define an incidence matrix $\bar{\mathbf{G}} \in \mathds{R}^{|\bar{V}| \times |E|}$ for $G$, where
\begin{align}\label{}
[\bar{\mathbf{G}}]_{jk}=\left\{
\begin{array}{ll}
  1 & \hbox{if $\bar{e}_k=(j,l) \in E$ and $j>l$} \\
  -1 & \hbox{if $\bar{e}_k=(j,l)\in E$ and $j<l$} \\
  0 & \hbox{otherwise}
\end{array}
\right..
\end{align}

\subsection{The Optimization Problem}

\begin{figure}
  \centering
  \includegraphics[width=8cm]{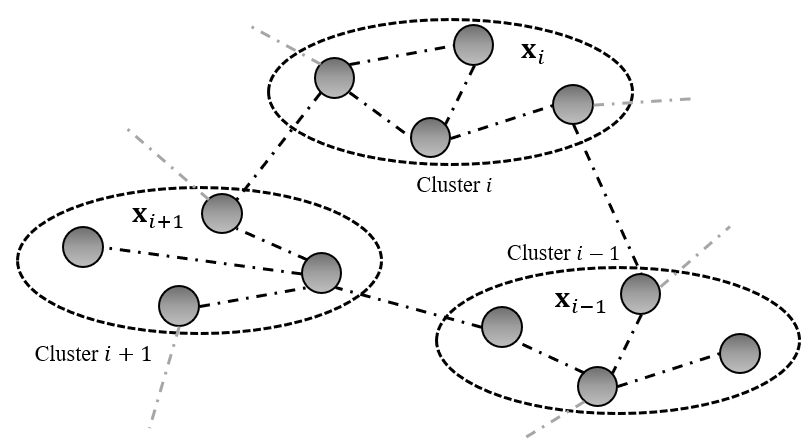}\\
  \caption{An illustrative communication topology of the network.}\label{clu}
\end{figure}

Let $H(\mathbf{x})= \sum_{i \in {V}} H_{i}(\mathbf{x}_{i})$ be the global cost function of network $G$ and $H_i(\mathbf{x}_i)= \sum_{j \in {V}_i} (f_{ij}(\mathbf{x}_{i}) + g_{ij}(\mathbf{x}_{i}))$ be the cost function of cluster $i$, where $f_{ij}+g_{ij}$ is the cost function of the $j$th agent in cluster $i$, $\mathbf{x}_{i} \in \mathds{R}^M$, $\mathbf{x}=[\mathbf{x}^{\top}_{1},...,\mathbf{x}^{\top}_{N}]^{\top} \in \mathds{R}^{{N}M}$.
The optimization problem of the whole network is formulated as
\begin{align}
  \mathrm{(P1)} \quad  \min\limits_{\mathbf{x}} \quad  & \sum_{i \in {V}} \sum_{j \in {V_i}} (f_{ij}(\mathbf{x}_{i}) + g_{ij} (\mathbf{x}_{i})) \nonumber \\
   \hbox{subject to} \quad & \mathbf{A}\mathbf{x} \leq \mathbf{b},
   \end{align}
with $\mathbf{A} \in \mathds{R}^{B \times {N}M}$, $\mathbf{b} \in \mathds{R}^{B}$. An illustrative communication topology of the network is shown in Fig. \ref{clu}.

\begin{Remark}\label{r0}
The composite cost function $f_{ij} + g_{ij}$ is a generalization of many cost functions in practical problems, such as resource allocation problems \cite{beck2014fast}, regularization problems \cite{hans2009bayesian}, and support vector machines \cite{zhao2017scope}. For instance, $g_{ij}$ can be a regularization penalty term, such as $l_1$- and $l_2$-norms. Alternatively, we can consider a local feasible region $X_{ij} \subseteq \mathds{R}^M$ for the $j$th agent in cluster $i$, where $X_{ij}$ is non-empty, convex, and closed.
Then, we can let $g_{ij}$ be an indicator function $\mathds{I}_{X_{ij}}$, where $\mathds{I}_{X_{ij}} (\mathbf{x}_{i})= \left\{\begin{array}{ll}
                    0 & \hbox{if $\mathbf{x}_{i} \in X_{ij}$} \\
                    +\infty  & \hbox{otherwise}
                  \end{array}
                  \right.$
 \cite{notarnicola2016asynchronous}.
\end{Remark}

\begin{Remark}\label{r0-1}
For the comparison purpose, we consider a conventional composite DOP
\begin{align}
  \mathrm{(P1+)} \quad  \min\limits_{\mathbf{x}} \quad  & \sum_{i \in {V}} (f_{i}(\mathbf{x}_{i}) + g_{i} (\mathbf{x}_{i})) \nonumber \\
   \hbox{subject to} \quad & \mathbf{A}\mathbf{x} \leq \mathbf{b},
   \end{align}
where $f_i+g_i$ is the cost function of agent $i$. Then, compared with Problem (P1+), the new features of the multi-cluster based Problem (P1) are explained as follows.
\begin{itemize}
  \item Problem (P1) can be viewed as a generalization of Problem (P1+) by expanding the agent who manages $\mathbf{x}_i$ into cluster $i$ with distributed agents. In case there is only one agent in each cluster, Problem (P1) is equivalent to Problem (P1+).
  \item To optimize the global cost function, the agents in each cluster need to generate a consensual decision, e.g., $\mathbf{x}_i$, which is involved by the coupling constraint $\mathbf{A}\mathbf{x} \leq \mathbf{b}$. Therefore, in the distributed setup, the interactions among the agents within each cluster and across different clusters should be considered simultaneously.
\end{itemize}

\end{Remark}

\begin{Assumption}\label{a0}
The edges in ${G}$ and $G_i$ are undirected; Any two distinct vertices in ${G}$ and $G_i$ are linked by at least one path, $\forall i\in V$.
\end{Assumption}

\begin{Assumption}\label{a1}
Both $f_{ij}: \mathds{R}^M \rightarrow (-\infty,+\infty]$ and $g_{ij}: \mathds{R}^M \rightarrow (-\infty,+\infty]$ are proper, convex, and closed extended real-valued functions; $f_{ij}$ is $\sigma_{ij}$-strongly convex, $\sigma_{ij}>0$, $\forall i\in {V}, j\in V_i$.
\end{Assumption}
The assumptions in Assumption \ref{a1} are often discussed in composite optimization problems \cite{shi2015proximal,wang2021distributed,schmidt2011convergence,chang2014multi,florea2020generalized,beck2014fast,notarnicola2016asynchronous}.

\begin{Assumption}\label{a3a}
{{(Constraint Qualification \cite{boyd2004convex})}} There exists an $\breve{\mathbf{x}} \in \mathbf{relint}D$ such that $\mathbf{A}\breve{\mathbf{x}} \leq \mathbf{b}$, where $D= \prod_{i\in V} \bigcap_{j \in V_i} D_{ij}$ with $D_{ij}$ the domain of $f_{ij} + g_{ij}$, $\forall i\in V$, $j\in V_i$.
\end{Assumption}

\section{Main Result}\label{sa4}

In this section, we propose a CDPG algorithm for solving the problem of interest, discuss the computational complexity with simple-structured cost functions, and further conduct the convergence analysis.

\subsection{Dual Problem}\label{sec1}

To realize distributed computation, we decouple the variable of clusters by defining $\mathbf{y}_{ij} \in \mathds{R}^M$ as the estimate of $\mathbf{x}_i$ by the $j$th agent in cluster $i$. Then, the collection of the estimates in cluster $i$ can be $\mathbf{y}_i=[\mathbf{y}^{\top}_{i1},...,\mathbf{y}^{\top}_{in_i}]^{\top} \in \mathds{R}^{n_iM}$ and the collection of the overall estimates can be $\mathbf{y} = [\mathbf{y}^{\top}_1,...,\mathbf{y}^{\top}_N]^{\top} \in \mathds{R}^{\sum_{i\in V}n_iM}$.

Based on Problem (P1), the estimates in each cluster should reach a consensus. Then by the consensus protocol in cluster $i$: $\mathbf{x}_i = \mathbf{y}_{i1} = ... = \mathbf{y}_{in_i}$, the global constraint can be represented by
\begin{align}\label{ab1}
\mathbf{A}\mathbf{x} \leq \mathbf{b}
& \Leftrightarrow  \sum_{i \in V} \mathbf{A}_i {\mathbf{x}}_i \leq  \mathbf{b} \Leftrightarrow  \sum_{i \in V} \frac{\mathbf{1}^{\top}_{n_i} \otimes \mathbf{A}_i}{n_i} \mathbf{y}_i \leq  \mathbf{b} \nonumber \\
& \Leftrightarrow  \sum_{i \in V} {\mathsf{A}_i\mathbf{y}_i} \leq  \mathbf{b} \Leftrightarrow  {\mathsf{A}\mathbf{y}} \leq  \mathbf{b},
\end{align}
where $\mathsf{A}_i=\frac{\mathbf{1}^{\top}_{n_i} \otimes \mathbf{A}_i}{n_i} \in \mathds{R}^{B \times n_iM}$ and $\mathsf{A} = [{\mathsf{A}_1},...,{\mathsf{A}_N}] \in \mathds{R}^{B \times \sum_{i \in V}n_i M}$ with $\mathbf{A}_i \in \mathds{R}^{B \times M}$ the $i$th column block of $\mathbf{A}$ (i.e., $\mathbf{A}=[\mathbf{A}_1,...,\mathbf{A}_i,...,\mathbf{A}_N]$). Essentially, (\ref{ab1}) reconstructs the coupling constraint in (P1) with the augmented variables without affecting the nature of the constraint when certain consensus is achieved in each cluster.

Note that the consensus constraint of $\mathbf{y}_{ij}$ in cluster $i$ can be equivalently written as $\mathsf{L}^i\mathbf{y}_i=\mathbf{0}$, where $\mathsf{L}^i = \mathbf{L}^i \otimes \mathbf{I}_M \in \mathds{R}^{n_iM \times n_iM}$ is the augmented Laplacian matrix of $G_i$ \cite{dimarogonas2010stability}. Then the consensus-based optimization problem of the whole network can be formulated as
\begin{align}
  \mathrm{(P2)} \quad  \min\limits_{\mathbf{y}} \quad  & \sum_{i \in {V}} \sum_{j \in {V}_i} (f_{ij}(\mathbf{y}_{ij}) + g_{ij} (\mathbf{z}_{ij})) \nonumber \\
   \hbox{subject to} \quad & \mathbf{y}_{ij} = \mathbf{z}_{ij}, \quad \forall i \in V, j \in V_i,  \\
   & \mathsf{L}^i\mathbf{y}_i=\mathbf{0},\quad \forall i \in V,\\
   & \mathsf{A} \mathbf{y} \leq \mathbf{b},
   \end{align}
where $\mathbf{z}_{ij} \in \mathds{R}^{M}$ is a slack variable.

The Lagrangian function of Problem (P2) can be given by
\begin{align}\label{11}
& L(\mathbf{y}, \mathbf{z}, \bm{\mu}, \bm{\nu},\bm{\phi})\nonumber \\
= & \sum_{i \in {V}} (\sum_{j \in {V}_i} (f_{ij}(\mathbf{y}_{ij}) + g_{ij}(\mathbf{z}_{ij}) +  \bm{\mu}_{ij}^{\top} (\mathbf{y}_{ij} - \mathbf{z}_{ij}))   \nonumber \\
& + \bm{\nu}^{\top}_{i}\mathsf{L}^i \mathbf{y}_i )+ \bm{\phi}^{\top} (\mathsf{A}\mathbf{y} - \mathbf{b}) \nonumber \\
= & \sum_{i \in {V}} (\sum_{j \in {V}_i} (f_{ij}(\mathbf{y}_{ij}) + g_{ij}(\mathbf{z}_{ij}) + \bm{\mu}_{ij}^{\top} (\mathbf{y}_{ij} - \mathbf{z}_{ij})) \nonumber \\
&  + \bm{\nu}^{\top}_{i}\mathsf{L}^i \mathbf{y}_i + \bm{\phi}^{\top} \mathsf{A}_i\mathbf{y}_i - \kappa_i \bm{\phi}^{\top}  \mathbf{b}) \nonumber \\
= & \sum_{i \in {V}} \sum_{j \in {V}_i} (f_{ij}(\mathbf{y}_{ij}) + g_{ij}(\mathbf{z}_{ij}) + \bm{\mu}_{ij}^{\top} (\mathbf{y}_{ij} - \mathbf{z}_{ij} ) \nonumber \\
& + \bm{\nu}^{\top}_{i}\mathsf{L}^i_{j} \mathbf{y}_{ij} + \bm{\phi}^{\top} \mathsf{A}_{ij}\mathbf{y}_{ij}  - \kappa_i \eta_{ij} \bm{\phi}^{\top}\mathbf{b}),
\end{align}
where $\bm{\mu}_{ij} \in \mathds{R}^M$ and $\bm{\nu}_i \in \mathds{R}^{n_iM}$ and $\bm{\phi} \in \mathds{R}^{B}$ are Lagrangian multipliers,
\begin{align}\label{}
&\mathbf{z}_{i} = [\mathbf{z}^{\top}_{i1},..., \mathbf{z}^{\top}_{in_i}]^{\top},\mathbf{z} = [\mathbf{z}^{\top}_{1},...,\mathbf{z}^{\top}_{N}]^{\top}, \nonumber \\
&\bm{\mu}_{i} = [\bm{\mu}^{\top}_{i1},..., \bm{\mu}^{\top}_{in_i}]^{\top},\bm{\mu} = [\bm{\mu}^{\top}_{1},...,\bm{\mu}^{\top}_{N}]^{\top}, \nonumber \\
&\bm{\nu} = [\bm{\nu}^{\top}_{1},...,\bm{\nu}^{\top}_{N}]^{\top}, \sum_{i \in {V}}\kappa_i = 1,\sum_{j \in {V}_i} \eta_{ij} = 1, \nonumber
\end{align}
$\mathsf{A}_{ij} \in \mathds{R}^{B \times M}$ is the $j$th column block of $\mathsf{A}_{i}$ (i.e., $\mathsf{A}_{i}=[\mathsf{A}_{i1},...,\mathsf{A}_{ij},...,\mathsf{A}_{in_i}]$), and $\mathsf{L}^i_{j} \in \mathds{R}^{n_iM \times M}$ is the $j$th column block of $\mathsf{L}^i$ (i.e., $\mathsf{L}^i=[\mathsf{L}^i_{1},...,\mathsf{L}^i_{j},...,\mathsf{L}^i_{n_i}]$).

Then the dual function can be obtained by
\begin{align}\label{}
& W(\bm{\mu},\bm{\nu},\bm{\phi}) \nonumber \\
= & \min\limits_{\mathbf{y}, \mathbf{z}} \sum_{i \in {V}} \sum_{j \in {V}_i} (f_{ij}(\mathbf{y}_{ij}) + g_{ij}(\mathbf{z}_{ij})  + \bm{\mu}_{ij}^{\top} (\mathbf{y}_{ij} - \mathbf{z}_{ij} ) \nonumber \\
&+ \bm{\nu}^{\top}_{i}\mathsf{L}^i_{j} \mathbf{y}_{ij} + \bm{\phi}^{\top} \mathsf{A}_{ij}\mathbf{y}_{ij}  - \kappa_i \eta_{ij} \bm{\phi}^{\top} \mathbf{b} )  \nonumber \\
= & \min\limits_{\mathbf{y}, \mathbf{z}} \sum_{i \in {V}} \sum_{j \in {V}_i} (f_{ij}(\mathbf{y}_{ij}) + (\bm{\mu}_{ij} + \mathsf{L}^{i\top}_{j}\bm{\nu}_{i} + \mathsf{A}^{\top}_{ij} \bm{\phi})^{\top} \mathbf{y}_{ij} \nonumber \\
& - \kappa_i \eta_{ij} \mathbf{b}^{\top} \bm{\phi} + g_{ij}(\mathbf{z}_{ij}) - \bm{\mu}^{\top}_{ij}\mathbf{z}_{ij})  \nonumber \\
= & \sum_{i \in {V}} \sum_{j \in {V}_i} (- f^{\diamond}_{ij}(-\bm{\mu}_{ij} - \mathsf{L}^{i\top}_{j}\bm{\nu}_{i} - \mathsf{A}^{\top}_{ij} \bm{\phi}) \nonumber \\
& - \kappa_i \eta_{ij} \mathbf{b}^{\top} \bm{\phi} - g_{ij}^{\diamond}(\bm{\mu}_{ij}))  \nonumber \\
= & \sum_{i \in {V}} \sum_{j \in {V}_i} (- f^{\diamond}_{ij}(\mathbf{H}_{ij} \bm{\vartheta}_{ij}) -\mathbf{E}_{ij} \bm{\vartheta}_{ij} - g_{ij}^{\diamond}(\mathbf{F}_{ij} \bm{\vartheta}_{ij})),
\end{align}
where
\begin{align}\label{}
  & \mathbf{H}_{ij}=[-\mathbf{I}_M,-\mathsf{L}^{i\top}_{j}, -\mathsf{A}^{\top}_{ij}],  \mathbf{E}_{ij}= [\mathbf{0}^{\top}_M,\mathbf{0}^{\top}_{n_iM}, \kappa_i \eta_{ij}\mathbf{b}^{\top}], \nonumber \\
  & \mathbf{F}_{ij}=[\mathbf{I}_M,\mathbf{O}_{M \times n_iM}, \mathbf{O}_{M \times B}], \bm{\vartheta}_{ij}=[\bm{\mu}_{ij}^{\top},\bm{\nu}^{\top}_{i}, \bm{\phi}^{\top}]^{\top}. \nonumber
\end{align}
Hence, the dual problem of Problem (P2) is
\begin{align}
\mathrm{(P3)} \quad & \min\limits_{\bm{\vartheta}, \bm{\phi} \geq \mathbf{0}} \sum_{i \in {V}} \sum_{j \in {V}_i} (f^{\diamond}_{ij}(\mathbf{H}_{ij} \bm{\vartheta}_{ij}) +\mathbf{E}_{ij} \bm{\vartheta}_{ij} + g_{ij}^{\diamond}(\mathbf{F}_{ij} \bm{\vartheta}_{ij})), \nonumber
\end{align}
where $\bm{\vartheta} = [\bm{\vartheta}^{\top}_{1},...,\bm{\vartheta}^{\top}_{N}]^{\top},\bm{\vartheta}_{i} = [\bm{\vartheta}^{\top}_{i1},..., \bm{\vartheta}^{\top}_{in_i}]^{\top}$. Since considering the constraint $\bm{\phi} \geq \mathbf{0}$ is equivalent to accommodating indictor function $\mathds{I}_{\mathds{R}^B_+}(\bm{\phi}) = \mathds{I}_{\mathds{R}^B_+}(\mathbf{N}_{ij}\bm{\vartheta}_{ij})$ into the cost function with $\mathbf{N}_{ij} = [\mathbf{O}_{B \times (M+n_iM)},\mathbf{I}_B]$, then Problem (P3) can be rewritten as
\begin{align}
\mathrm{(P4)} \quad \min\limits_{\bm{\vartheta}} & \sum_{i \in {V}} \sum_{j \in {V}_i} p_{ij}(\bm{\vartheta}_{ij}) + \sum_{i \in {V}} \sum_{j \in {V}_i} q_{ij}(\bm{\vartheta}_{ij}), \nonumber
\end{align}
where
\begin{align}\label{}
& p_{ij}(\bm{\vartheta}_{ij}) = f^{\diamond}_{ij}(\mathbf{H}_{ij} \bm{\vartheta}_{ij}) +\mathbf{E}_{ij} \bm{\vartheta}_{ij}, \\
& q_{ij}(\bm{\vartheta}_{ij}) = g_{ij}^{\diamond}(\mathbf{F}_{ij} \bm{\vartheta}_{ij}) + \mathds{I}_{\mathds{R}^B_+}(\mathbf{N}_{ij}\bm{\vartheta}_{ij}).
\end{align}

Note that $\bm{\vartheta}_{ij}$ is coupled among the agents with the common components $\bm{\nu}_i$ and $\bm{\phi}$. To construct a separable structure for the cost function in (P4), we define $\bm{\lambda}_{ij}=[\bm{\mu}_{ij}^{\top},\bm{\gamma}^{\top}_{ij},\bm{\theta}_{ij}^{\top}]^{\top}$, where $\bm{\gamma}_{ij}$ and $\bm{\theta}_{ij}$ are the local estimates of $\bm{\nu}_i$ and $\bm{\phi}$ by the $j$th agent in cluster $i$, respectively.
In addition, we define
\begin{align}\label{}
& \bm{\gamma}_{i}=[\bm{\gamma}^{\top}_{i1},..., \bm{\gamma}^{\top}_{in_i}]^{\top}, \bm{\gamma}=[\bm{\gamma}^{\top}_{1},..., \bm{\gamma}^{\top}_{N}]^{\top}, \nonumber \\
& \bm{\theta}_{i}= [\bm{\theta}^{\top}_{i1},..., \bm{\theta}^{\top}_{in_i}]^{\top},  \bm{\theta}=[\bm{\theta}^{\top}_{1},..., \bm{\theta}^{\top}_{N}]^{\top}, \nonumber \\
& \bm{\lambda}_i=[\bm{\lambda}^{\top}_{i1} ,...,\bm{\lambda}^{\top}_{in_i}]^{\top}, \bm{\lambda}=[\bm{\lambda}^{\top}_1 ,...,\bm{\lambda}^{\top}_{N}]^{\top}, \nonumber \\
& \mathbf{M}_{ij} = [\mathbf{O}_{n_iM \times M},\mathbf{I}_{n_iM},\mathbf{O}_{n_iM \times B}], \nonumber \\
& \mathbf{M}_i= \mathbf{I}_{n_i} \otimes \mathbf{M}_{ij}, \mathbf{M}=\mathrm{diag}[\mathbf{M}_{i}]_{V}, \nonumber \\
&\mathbf{N}_i = \mathbf{I}_{n_i} \otimes \mathbf{N}_{ij}, \mathbf{N}=\mathrm{diag}[\mathbf{N}_{i}]_{V}.\nonumber
\end{align}
Then we have $\bm{\gamma}_{ij} = \mathbf{M}_{ij}\bm{\lambda}_{ij}$, $\bm{\theta}_{ij} =  \mathbf{N}_{ij}\bm{\lambda}_{ij}$, $\bm{\gamma}_i = \mathbf{M}_i\bm{\lambda}_i$, $\bm{\gamma} = \mathbf{M}\bm{\lambda}$, $\bm{\theta}_i =  \mathbf{N}_i\bm{\lambda}_i$, and $\bm{\theta} =  \mathbf{N}\bm{\lambda}$.
For convenience purposes, we relabel $\bm{\theta}_{ij}$ by defining $\bar{\bm{\theta}}_{n_{ij}}=\bm{\theta}_{ij}$. Then the consensus of ${\bm{\gamma}}_{ij}$ in $V_i$ and $\bar{\bm{\theta}}_{n_{ij}}$ in $\bar{V}$ can be characterized by
\begin{align}\label{}
\bm{\gamma}_{il} & - \bm{\gamma}_{ij} =\mathbf{0}, \quad \forall l \in S_{ij}, \label{a2} \\
\bar{\bm{\theta}}_{k} & - \bar{\bm{\theta}}_{n_{ij}} = \mathbf{0}, \quad  \forall k \in \bar{S}_{ij}, \label{a3}
\end{align}
respectively, where $S_{ij} = \{l| (j,l) \in E_i, l >j\}$, $\bar{S}_{ij} = \{k| (n_{ij},k) \in E, k> n_{ij} \}$, $\forall i \in V$, $j \in V_i$.

\begin{Example}\label{e31}
Assume that the cluster setting in Example \ref{e1} follows $\bar{V} = \{1,2,3,4,5\}$, $V_1=\{1\}$, $V_2=\{1\}$, $V_3=\{1\}$, $V_4=\{1,2\}$. Then $S_{11}=\emptyset$, $S_{21}=\emptyset$, $S_{31}=\emptyset$, $S_{41}=\{2\}$, and $S_{42}=\emptyset$. In addition, $\bar{S}_{11}=\{2\}$, $\bar{S}_{21}=\{4,5\}$, $\bar{S}_{31}=\{4\}$, $\bar{S}_{41}=\{5\}$, and $\bar{S}_{42}=\emptyset$. An empty set of $S_{ij}$ or $\bar{S}_{ij}$ means no consensus constraint exists in (\ref{a2}) or (\ref{a3}), respectively.
\end{Example}

Note that (\ref{a2}) and (\ref{a3}) can be written in compact forms with the help of the incidence matrix \cite{dimarogonas2010stability}. By defining $\mathsf{G}_i = \mathbf{G}_i^{\top} \otimes \mathbf{I}_{n_iM}$, (\ref{a2}) can be represented by $\mathsf{G}_i \bm{\gamma}_{i}= \mathsf{G}_i \mathbf{M}_i \bm{\lambda}_{i} = \mathbf{0}$. In addition, one can construct $\mathsf{G}\bm{\gamma} = \mathsf{G} \mathbf{M} \bm{\lambda} = \mathbf{0}$ by including all clusters with $\mathsf{G}=\mathrm{diag}[\mathsf{G}_i]_{V}$. Similarly, (\ref{a3}) can be represented by $\bar{\mathsf{G}} \bar{\bm{\theta}}= \bar{\mathsf{G}}\mathbf{N}{\bm{\lambda}} = \mathbf{0}$, where $\bar{\bm{\theta}}=[\bar{\bm{\theta}}^{\top}_1,...,\bar{\bm{\theta}}^{\top}_{\sum_{l=1}^Nn_l}]^{\top}$ and $\bar{\mathsf{G}} = \bar{\mathbf{G}}^{\top} \otimes \mathbf{I}_B$. Then, by defining $\mathbf{Z} = [ \mathbf{M}^{\top}\mathsf{G}^{\top},\mathbf{N}^{\top}\bar{\mathsf{G}}^{\top}]^{\top}$,
(\ref{a2}) and (\ref{a3}) can be jointly represented by $\mathbf{Z}\bm{\lambda}=\mathbf{0}$.

Consider a quadratic term
\begin{align}\label{}
\|\mathbf{Z}\bm{\lambda} \|_{\mathbf{D}[\pi]}^2
= & \sum_{i \in V} \sum_{j \in V_i} \pi_{ij} \sum_{ l \in S_{ij}} (\bm{\gamma}_{il} - \bm{\gamma}_{ij})^2 \nonumber \\
 & + \sum_{i \in V} \sum_{j \in V_i}  \pi_{ij} \sum_{ k \in \bar{S}_{ij}} (\bar{\bm{\theta}}_{k}- \bar{\bm{\theta}}_{n_{ij}})^2, \nonumber
\end{align}
where
\begin{align}\label{pi}
& \mathbf{D}[\pi] = \left[
                           \begin{array}{c|c}
                             \mathrm{diag}[\mathbf{V}_i]_V & \mathbf{O} \\
                             \hline
                              \mathbf{O} & \mathrm{diag}[\bar{\mathbf{V}}_i]_V \\
                           \end{array}
                         \right], \nonumber \\
& \mathbf{V}_i=  \left[
                      \begin{array}{cccc}
                        \pi_{i1}\mathbf{I}_{|S_{i1}|n_iM} & & & \mathbf{O} \\
                         & \pi_{i2}\mathbf{I}_{|S_{i2}|n_iM}  & & \\
                         & & \ddots & \\
                         \mathbf{O}& & & \pi_{in_i}\mathbf{I}_{|S_{in_i}|n_iM}  \\
                      \end{array}
                    \right], \nonumber \\
&  \bar{\mathbf{V}}_i= \left[
                      \begin{array}{cccc}
                        \pi_{i1}\mathbf{I}_{|\bar{S}_{i1}|B} & & & \mathbf{O} \\
                         &  \pi_{i2}\mathbf{I}_{|\bar{S}_{i2}|B} & & \\
                         &  & \ddots & \\
                         \mathbf{O}&  & & \pi_{in_i}\mathbf{I}_{|\bar{S}_{in_i}|B}  \\
                      \end{array}
                    \right], \nonumber
\end{align}
with $\pi_{ij}>0$ a weighting factor, $i\in V$, $j\in V_i$.

Then, a consensus-based penalized dual problem of Problem (P2) can be given by
\begin{align}
\mathrm{(P5)} \quad \min\limits_{ {\bm{\lambda}} } \quad & \Phi ( {\bm{\lambda}}) + \frac{1}{2} \|\mathbf{Z}\bm{\lambda} \|_{\mathbf{D}[\pi]}^2 \nonumber \\
\hbox{subject to} \quad & \mathbf{Z}\bm{\lambda} = \mathbf{0},
\end{align}
where
\begin{align}\label{}
\Phi ( {\bm{\lambda}}) & = P( {\bm{\lambda}}) + Q ( {\bm{\lambda}}), \\
P( {\bm{\lambda}}) & = \sum_{i \in {V}} \sum_{j \in {V}_i} p_{ij}(\bm{\lambda}_{ij}), \\
Q( {\bm{\lambda}}) & = \sum_{i \in {V}} \sum_{j \in {V}_i} q_{ij}(\bm{\lambda}_{ij}), \\
p_{ij}(\bm{\lambda}_{ij}) & = f^{\diamond}_{ij}(\mathbf{H}_{ij} \bm{\lambda}_{ij}) +\mathbf{E}_{ij} \bm{\lambda}_{ij}, \\
q_{ij}(\bm{\lambda}_{ij}) & = g_{ij}^{\diamond}(\mathbf{F}_{ij} \bm{\lambda}_{ij}) + \mathds{I}_{\mathds{R}^B_+}(\mathbf{N}_{ij}\bm{\lambda}_{ij}). \label{21}
\end{align}
$\Phi$ is convex based on the definition of Fenchel conjugate.

\subsection{Distributed Optimization Algorithm Development}\label{sec1}

The Lagrangian function of Problem (P5) can be given by
 \begin{align}\label{fp1}
\mathcal{L} (\bm{\lambda},  \bm{\omega}) = & P(\bm{\lambda})+ Q(\bm{\lambda}) + \frac{1}{2}\|\mathbf{Z}\bm{\lambda}\|_{\mathbf{D}[\pi]}^2 + \bm{\omega}^{\top}\mathbf{Z} \bm{\lambda} \nonumber \\
= & \sum_{i \in {V}} \sum_{j \in {V}_i} p_{ij}(\bm{\lambda}_{ij}) + \sum_{i \in {V}} \sum_{j \in {V}_i} q_{ij}(\bm{\lambda}_{ij}) \nonumber \\
 & + \frac{1}{2} \sum_{i \in V} \sum_{j \in V_i} \pi_{ij} \sum_{l \in S_{ij}} (\bm{\gamma}_{il} - \bm{\gamma}_{ij})^2 \nonumber \\
 & + \frac{1}{2} \sum_{i \in V} \sum_{j \in V_i}  \pi_{ij} \sum_{ k \in \bar{S}_{ij}} (\bar{\bm{\theta}}_{k}- \bar{\bm{\theta}}_{n_{ij}})^2 \nonumber \\
& + \sum_{i\in V} \sum_{j \in V_i} \sum_{l \in S_{ij}} \bm{\xi}_{ijl}^{\top} (\bm{\gamma}_{il} - \bm{\gamma}_{ij}) \nonumber \\
 & + \sum_{i\in V} \sum_{j \in V_i} \sum_{k \in \bar{S}_{ij}} \bm{\zeta}_{ijk}^{\top} (\bar{\bm{\theta}}_{k}- \bar{\bm{\theta}}_{n_{ij}}),
\end{align}
where
\begin{align}\label{}
& \bm{\xi}_{ijl} \in \mathds{R}^{n_iM}, \bm{\xi}_{ij} = [ \bm{\xi}^{\top}_{ijl_1},...,\bm{\xi}^{\top}_{ijl_{|S_{ij}|}}]^{\top}, l_{(\cdot)} \in S_{ij}, \nonumber\\
& \bm{\zeta}_{ijk} \in \mathds{R}^B, \bm{\zeta}_{ij} =[\bm{\zeta}^{\top}_{ijk_1},...,\bm{\zeta}^{\top}_{ijk_{|\bar{S}_{ij}|}}]^{\top}, k_{(\cdot)} \in \bar{S}_{ij}, \nonumber \\
& \bm{\xi}_{i} = [\bm{\xi}^{\top}_{i1},...,\bm{\xi}^{\top}_{in_i}]^{\top}, \bm{\xi} = [\bm{\xi}^{\top}_{1},...,\bm{\xi}^{\top}_{N}]^{\top}, \nonumber \\
& \bm{\zeta}_{i} = [\bm{\zeta}^{\top}_{i1},...,\bm{\zeta}^{\top}_{in_i}]^{\top}, \bm{\zeta} = [\bm{\zeta}^{\top}_{1},...,\bm{\zeta}^{\top}_{N}]^{\top}, \bm{\omega}=[\bm{\xi}^{\top},\bm{\zeta}^{\top}]^{\top}. \nonumber
\end{align}
Here, $\bm{\omega}$ is the collection of the Lagrangian multipliers. Sequence $\{l_1,...,l_{|S_{ij}|}\}$ is decided by: $\forall l_m,l_n \in S_{ij}$, if $l_m>l_n$, then $m>n$, and vice versa. Similarly, in $\{k_1,...,k_{|\bar{S}_{ij}|}\}$, $\forall k_m,k_n \in \bar{S}_{ij}$, if $k_m>k_n$, then $m>n$, and vice versa.

\begin{Example}\label{e4}
Based on the graph in Example \ref{e1} and the cluster setting in Example \ref{e31}, we have $\bm{\xi}_{41}=\bm{\xi}_{412}$, $\bm{\zeta}_{11}=\bm{\zeta}_{112}$, $\bm{\zeta}_{21}=[\bm{\zeta}^{\top}_{214},\bm{\zeta}^{\top}_{215}]^{\top}$, $\bm{\zeta}_{31}=\bm{\zeta}_{314}$, and $\bm{\zeta}_{41}=\bm{\zeta}_{415}$.
\end{Example}

Let ${C}$ be the set of the saddle points of $\mathcal{L} (\bm{\lambda},\bm{\omega})$. Then, any saddle point $(\bm{\lambda}^*,\bm{\omega}^*) \in {C}$ satisfies \cite{rockafellar1970convex}
\begin{align}\label{sad}
 \mathcal{L} (\bm{\lambda},\bm{\omega}^*) \geq  \mathcal{L} (\bm{\lambda}^*,\bm{\omega}^*) \geq  \mathcal{L} (\bm{\lambda}^*,\bm{\omega}).
\end{align}
We aim to seek a saddle point of $\mathcal{L}(\bm{\lambda},\bm{\omega}) $, which can be characterized by Karush-Kuhn-Tucker (KKT) conditions  \cite{hanson1981sufficiency}
\begin{align}
 & \mathbf{0}\in \nabla_{\bm{\lambda}} P(\bm{\lambda}^*) + \partial_{\bm{\lambda}} Q(\bm{\lambda}^*)
 + \mathbf{Z}^{\top} \bm{\omega}^* + \mathbf{Z}^{\top}\mathbf{D}[\pi]\mathbf{Z}\bm{\lambda}^*,    \label{k1} \\
&  \mathbf{Z} \bm{\lambda}^* = \mathbf{0}. \label{k2}
\end{align}
Then the proposed CDPG algorithm for solving Problem (P5) is designed as
\begin{align}
 \bm{\lambda}^{t+1} = & \mathrm{prox}^{\mathbf{S}[c]}_{Q}[\bm{\lambda}^t - \mathbf{S}[c] ( \nabla_{\bm{\lambda}} P (\bm{\lambda}^t) + {\mathbf{Z}}^{\top} \bm{\omega}^t \nonumber \\
& \quad \quad \quad \quad  + \mathbf{Z}^{\top}\mathbf{D}[\pi]\mathbf{Z}\bm{\lambda}^t)], \label{f1} \\
 \bm{\omega}^{t+1}= & \bm{\omega}^t + \mathbf{D}[\pi] {\mathbf{Z}} \bm{\lambda}^{t+1}, \label{f2}
\end{align}
where
\begin{align}\label{cp}
& \mathbf{S}[c] = \mathrm{diag}[\mathbf{C}_i]_V, \\
& \mathbf{C}_i= \left[
                      \begin{array}{cccc}
                        c_{i1} &   & & \mathbf{O} \\
                        &   c_{i2} & & \\
                        &  & \ddots & \\
                       \mathbf{O}  &  &  & c_{in_i}  \\
                      \end{array}
                    \right] \otimes \mathbf{I}_{M+n_iM+B}.
\end{align}
Based on (\ref{21}), we have
\begin{align}\label{}
q_{ij} (\bm{\lambda}_{ij})& = g_{ij}^{\diamond}(\bm{\mu}_{ij}) + \mathds{I}_{\mathds{R}^B_+}(\bm{\theta}_{ij}) \nonumber \\
& = g_{ij}^{\diamond}(\bm{\mu}_{ij}) + \mathds{I}_{\mathds{R}^{n_iM}}(\bm{\gamma}_{ij}) + \mathds{I}_{\mathds{R}^B_+}(\bm{\theta}_{ij}).
\end{align}
Then one can have \cite[Thm. 6.6]{beck2017first}
\begin{align}\label{34}
\mathrm{prox}^{c_{ij}}_{q_{ij}} =& \mathrm{prox}^{c_{ij}}_{g_{ij}^{\diamond}} \times \mathrm{prox}^{c_{ij}}_{\mathds{I}_{\mathds{R}^{n_iM}}} \times \mathrm{prox}^{c_{ij}}_{\mathds{I}_{\mathds{R}^B_+}}.
\end{align}
Based on (\ref{34}), (\ref{f1}) and (\ref{f2}) can be decomposed into
\begin{align}
 \bm{\mu}_{ij}^{t+1}  = & \mathrm{prox}^{c_{ij}}_{g_{ij}^{\diamond}} [\bm{\mu}_{ij}^t - c_{ij} \nabla_{\bm{\mu}_{ij}} p_{ij} (\bm{\lambda}_{ij}^t) ], \label{1} \\
 \bm{\gamma}_{ij}^{t+1}  = & \bm{\gamma}_{ij}^t - c_{ij}  ( \nabla_{\bm{\gamma}_{ij}} p_{ij} (\bm{\lambda}_{ij}^t) -  \sum_{ l \in S_{ij}} \bm{\xi}_{ijl}^t  \nonumber \\
 &  + \sum_{ l^{'} \in S^{\sharp}_{ij}} \bm{\xi}_{il^{'}j}^t + \pi_{ij} \sum_{l \in S_{ij}}  ( \bm{\gamma}_{ij}^t - \bm{\gamma}_{il}^t)  \nonumber \\
 & + \sum_{l^{'} \in S^{\sharp}_{ij}} \pi_{il^{'}}  (\bm{\gamma}_{ij}^t - \bm{\gamma}_{il^{'}}^t )) ,\label{2} \\
 \bm{\theta}_{ij}^{t+1} = & \mathrm{P}_{\mathbb{R}^B_+} [\bm{\theta}_{ij}^t - c_{ij}  ( \nabla_{\bm{\theta}_{ij}} p_{ij} (\bm{\lambda}_{ij}^t)  - \sum_{k \in \bar{S}_{ij}}
{\bm{\zeta}}_{ijk}^t \nonumber \\
 & +  \sum_{ k^{'} \in \bar{S}^{\sharp}_{ij}} {\bm{\zeta}}_{ik^{'}j}^t  + \pi_{ij}\sum_{k \in \bar{S}_{ij} } (\bar{\bm{\theta}}_{n_{ij}}^t - \bar{\bm{\theta}}_{k}^t  )  \nonumber \\
 &  + \sum_{k^{'} \in \bar{S}^{\sharp}_{ij}} \bar{\pi}_{k^{'}} (\bar{\bm{\theta}}_{n_{ij}}^t- \bar{\bm{\theta}}_{k^{'}}^t))], \label{3} \\
 \bm{\xi}_{ijl}^{t+1} = &  \bm{\xi}_{ijl}^t  + \pi_{ij} (\bm{\gamma}_{il}^{t+1} - \bm{\gamma}_{ij}^{t+1}), \quad \forall l \in {S}_{ij} , \label{4} \\
  \bm{\zeta}_{ijk}^{t+1} =  & \bm{\zeta}_{ijk}^t + \pi_{ij} (\bar{\bm{\theta}}_{k}^{t+1}- \bar{\bm{\theta}}_{n_{ij}}^{t+1}), \quad \forall k \in \bar{S}_{ij}, \label{5}
\end{align}
due to the separability of $P$ and $Q$, where $\bar{\pi}$ is the relabeled $\pi$ by formula $\bar{\pi}_{n_{ij}} = \pi_{ij}$, $S^{\sharp}_{ij} = \{l| (j,l) \in E_i, l < j\}$, $\bar{S}^{\sharp}_{ij} = \{k| (n_{ij},k) \in E, k < n_{ij} \}$, and $\mathrm{P}_{\mathbb{R}^B_+}[\cdot]$ is a Euclidean projection onto $\mathbb{R}^B_+$, $\forall i \in {V}, j \in V_i$. (\ref{2}) and (\ref{3}) hold since the proximal mapping of an indicator function is equivalent to a Euclidean projection (the projection onto $\mathds{R}^{n_iM}$ in (\ref{2}) is omitted) \cite[Sec. 1.2]{parikh2014proximal}. 

\begin{Remark}\label{r3}
To apply (\ref{1}) to (\ref{3}), one needs to compute $\nabla_{\bm{\mu}_{ij}}p_{ij} (\bm{\lambda}_{ij})$, $\nabla_{\bm{\gamma}_{ij}}p_{ij} (\bm{\lambda}_{ij})$, and $\nabla_{\bm{\theta}_{ij}}p_{ij} (\bm{\lambda}_{ij})$. By Lemma \ref{l1}, we have
\begin{align}\label{}
\nabla_{\bm{\lambda}_{ij}} p_{ij} & (\bm{\lambda}_{ij}) = [\nabla^{\top}_{\bm{\mu}_{ij}}p_{ij} (\bm{\lambda}_{ij}),        \nabla^{\top}_{\bm{\gamma}_{ij}}p_{ij} (\bm{\lambda}_{ij}),        \nabla^{\top}_{\bm{\theta}_{ij}}p_{ij} (\bm{\lambda}_{ij})]^{\top} \nonumber \\ = & \mathbf{H}_{ij}^{\top}  \nabla_{\mathbf{H}_{ij} \bm{\lambda}_{ij}} f^{\diamond}_{ij}(\mathbf{H}_{ij} \bm{\lambda}_{ij}) + \mathbf{E}_{ij}^{\top} \nonumber \\
= & \mathbf{H}_{ij}^{\top} \arg \max\limits_{\mathbf{u}} ( (\mathbf{H}_{ij} \bm{\lambda}_{ij})^{\top} \mathbf{u} - f_{ij}(\mathbf{u}))  + \mathbf{E}_{ij}^{\top},  \label{6}
\end{align}
which can be completed by some numerical methods only with local information. On the other hand, if the proximal mapping of $g_{ij}$ is easier to obtain, the computation of (\ref{1}) can be further simplified by employing Lemma \ref{md}, which gives
\begin{align}
& \bm{\varrho}_{ij}^t  = \bm{\mu}_{ij}^t - c_{ij} \nabla_{\bm{\mu}_{ij}} p_{ij} (\bm{\lambda}_{ij}^t), \label{7} \\
& \bm{\mu}_{ij}^{t+1}   =  \mathrm{prox}^{c_{ij}}_{g_{ij}^{\diamond}} [ \bm{\varrho}_{ij}^t ] = \bm{\varrho}_{ij}^t - {c_{ij}}\mathrm{prox}^{\frac{1}{c_{ij}}}_{g_{ij}} [\frac{ \bm{\varrho}_{ij}^t}{c_{ij}} ], \label{8}
\end{align}
where we use $g^{\diamond \diamond }_{ij} = g_{ij} $ since $g_{ij}$ is convex and lower semi-continuous, and $g^{\diamond \diamond }_{ij}$ is the biconjugate of $g_{ij}$ \cite[Sec. 3.3.2]{boyd2004convex}. Specifically, if $g_{ij}=\mathds{I}_{X_{ij}}$ (see Remark \ref{r0}), (\ref{8}) is equivalent to
$\bm{\mu}_{ij}^{t+1} =  \bm{\varrho}_{ij}^t - {c_{ij}} \mathrm{P}_{X_{ij}}[\frac{ \bm{\varrho}_{ij}^t}{c_{ij}} ]$.
\end{Remark}

\begin{Remark}
The considered inequality-constrained DOP can cover the equality-constrained scenarios by removing the constraint $\bm{\phi} \geq \mathbf{0}$ in Problem (P3) based on KKT conditions, which leads to $q_{ij}(\bm{\lambda}_{ij}) = g_{ij}^{\diamond}(\mathbf{F}_{ij} \bm{\lambda}_{ij})$ in Problem (P5) by removing the indicator function $\mathds{I}_{\mathds{R}^B_+}(\mathbf{N}_{ij}\bm{\lambda}_{ij})$ from the cost function. Then, the updating law (\ref{3}) becomes
\begin{align}\label{}
 \bm{\theta}_{ij}^{t+1} = & \bm{\theta}_{ij}^t - c_{ij}  ( \nabla_{\bm{\theta}_{ij}} p_{ij} (\bm{\lambda}_{ij}^t)  - \sum_{k \in \bar{S}_{ij}}
{\bm{\zeta}}_{ijk}^t +  \sum_{ k^{'} \in \bar{S}^{\sharp}_{ij}} {\bm{\zeta}}_{ik^{'}j}^t  \nonumber \\
 &  + \pi_{ij}\sum_{k \in \bar{S}_{ij} } (\bar{\bm{\theta}}_{n_{ij}}^t - \bar{\bm{\theta}}_{k}^t  ) + \sum_{k^{'} \in \bar{S}^{\sharp}_{ij}} \bar{\pi}_{k^{'}} (\bar{\bm{\theta}}_{n_{ij}}^t- \bar{\bm{\theta}}_{k^{'}}^t)) \nonumber
\end{align}
by removing the Euclidean projection.
\end{Remark}

In the following, we discuss how to recover the optimal primal solution $\mathbf{y}^*$. By the saddle point property
\begin{align}\label{6}
L(\mathbf{y}^*, \mathbf{z}^*, \bm{\mu}, \bm{\nu},\bm{\phi}) & \leq  L(\mathbf{y}^*, \mathbf{z}^*, \bm{\mu}^*, \bm{\nu}^*,\bm{\phi}^*) \nonumber \\
& \leq L(\mathbf{y}, \mathbf{z}, \bm{\mu}^*, \bm{\nu}^*,\bm{\phi}^*),
\end{align}
the optimal primal variable can be obtained by the second inequality in (\ref{6}): $\mathbf{y}^* = \arg \min_{\mathbf{y}} L(\mathbf{y}, \mathbf{z}, \bm{\mu}^*, \bm{\nu}^*,\bm{\phi}^*)$. Then by decomposing $\mathbf{y}$ and omitting the constant terms in $L$, we have
\begin{align}\label{rr1}
\mathbf{y}_{ij}^* = & \arg \min_{\mathbf{y}_{ij}} f_{ij}(\mathbf{y}_{ij}) + (\bm{\mu}_{ij}^{*} + \mathsf{L}^{i\top}_{j} \bm{\nu}^{*}_{i} + \mathsf{A}^{\top}_{ij} \bm{\phi}^* )^{\top} \mathbf{y}_{ij} \nonumber \\
= & \arg \min_{\mathbf{y}_{ij}} f_{ij}(\mathbf{y}_{ij}) + (\bm{\mu}_{ij}^{*} + \mathsf{L}^{i\top}_{j} \bm{\gamma}_{ij}^{*} + \mathsf{A}^{\top}_{ij} \bm{\theta}_{ij}^{*})^{\top} \mathbf{y}_{ij},
\end{align}
where we use $\bm{\gamma}^*_{ij}= \bm{\nu}_i^{*}$ and $\bm{\theta}^*_{ij} = \bm{\phi}^*$ since $\bm{\gamma}^*_{ij}$ and $\bm{\theta}^*_{ij}$ are the optimal local estimates of $\bm{\nu}_i^{*}$ and $\bm{\phi}^*$, respectively, $i\in V$, $j \in V_i$. The detailed computation procedure of the CDPG algorithm is summarized in Algorithm \ref{ax1}.
\begin{algorithm}
\caption{CDPG algorithm}\label{ax1}
\begin{algorithmic}[1]
\State Initialize $\bm{\lambda}^0$, $\bm{\omega}^0$. Determine step-sizes $c_{ij},\pi_{ij}>0$, $\forall i \in {V}, j \in V_i$.
\For {$t= 0,1,2,...$}
\For {$i= 1,2,...,{N}$ and $j=1,2,...,n_i$} (in parallel)
\State Update $\nabla_{\bm{\mu}_{ij}}p_{ij} (\bm{\lambda}^t_{ij})$, $\nabla_{\bm{\gamma}_{ij}}p_{ij} (\bm{\lambda}^t_{ij})$, and $\nabla_{\bm{\theta}_{ij}}p_{ij} (\bm{\lambda}^t_{ij})$ based on (\ref{6}).
\State Update $\bm{\mu}^{t+1}_{ij}$ based on (\ref{1}) or (\ref{7})-(\ref{8}).
\State \textbf{Intra-cluster interaction}:
\State Update $\bm{\gamma}^{t+1}_{ij}$ and $\bm{\xi}^{t+1}_{ij}$ based on (\ref{2}) and (\ref{4}), respectively.
\State \textbf{Extra-cluster interaction}:
\State Update $\bm{\theta}^{t+1}_{ij}$ and $\bm{\zeta}^{t+1}_{ij}$ based on (\ref{3}) and (\ref{5}), respectively.
\EndFor
\EndFor
\State Obtain the output of dual variables ${\bm{\mu}}_{ij},$ ${\bm{\gamma}}_{ij}$, ${\bm{\xi}}_{ij}$, ${\bm{\theta}}_{ij}$, and ${\bm{\zeta}}_{ij}$, $\forall i \in V, j\in V_i$. Calculate the primal solution based on (\ref{rr1}).
\end{algorithmic}
\end{algorithm}


\subsection{Convergence Analysis}\label{sa5}

\begin{Lemma}\label{lam2}
With Assumption \ref{a1}, the Lipschitz constant of $\nabla_{\bm{\lambda}_{ij}}  p_{ij}(\bm{\lambda}_{ij})$ is given by $h_{ij}= \frac{ \| \mathbf{H}_{ij} \|^2 }{\sigma_{ij}}$, $\forall i \in V$, $j \in V_i$.
\end{Lemma}

See the proof in Appendix \ref{lam2p}.

In the following, we let $\bar{\mathbf{D}}[\pi]$ and $\bar{\mathbf{S}}[c]$ be the inverse matrices of ${\mathbf{D}}[\pi]$ and ${\mathbf{S}}[c]$, respectively. In addition, define
\begin{align}\label{cp}
& \mathbf{S}[h] = \mathrm{diag}[\mathbf{W}_i]_V, \\
& \mathbf{W}_i= \left[
                      \begin{array}{cccc}
                        h_{i1} &  & & \mathbf{O} \\
                        &  h_{i2} & & \\
                        & & \ddots  & \\
                       \mathbf{O}  & & & h_{in_i}  \\
                      \end{array}
                    \right] \otimes \mathbf{I}_{M+n_iM+B}.
\end{align}

\begin{Theorem}\label{th1}
Suppose that Assumptions \ref{a0}-\ref{a3a} hold. Let $0< c_{ij}  \leq \frac{1}{h_{ij} + \tau_{\mathrm{max}}(\mathbf{Z}^{\top}\mathbf{D}[\pi]\mathbf{Z})}$, $\forall i \in {V}, j\in V_i$. By Algorithm \ref{ax1}, for any $(\bm{\lambda}^*,\bm{\omega}^*) \in {C}$, we have
\begin{align}
&  | \Phi (\bar{\bm{\lambda}}^{T+1}) - \Phi(\bm{\lambda}^*) | \leq \frac{\Theta}{T+1} , \label{t1} \\
 & \| \bm{\omega}^*  \| \| \mathbf{Z} \bar{\bm{\lambda}}^{T+1} \|  \leq  \frac{\Theta}{ T+1} ,\label{t2}
\end{align}
where $\Theta = \|  \bm{\omega}^* \|^2_{4\bar{\mathbf{D}}[\pi]} + \| \bm{\omega}^0 \|^2_{\bar{\mathbf{D}}[\pi]}  + \|  \bm{\lambda}^* -\bm{\lambda}^0 \|^2_{\frac{1}{2}\bar{\mathbf{S}}[c] -\frac{1}{2}\mathbf{Z}^{\top}\mathbf{D}[\pi]\mathbf{Z}}$, $\bar{\bm{\lambda}}^{T+1} = \frac{1}{T+1}\sum_{t=0}^T {\bm{\lambda}}^{t+1}$, $T \in \mathds{N}_+$.
\end{Theorem}

See the proof in Appendix \ref{th1p}.


\subsection{Computational Complexity with Simple-Structured Cost Functions}\label{sa41}
By the proposed CDPG algorithm, (\ref{1}) to (\ref{3}) may require some inner-loop optimizations to compute $\nabla p_{ij}$ as discussed in Remark \ref{r3} and the proximal mapping of $g_{ij}^{\diamond}$. In the following, we will discuss the CDPG algorithm with some simple-structured cost functions, where the computational complexity can be reduced. For the gradient-based iterative algorithms, the computational complexity is dominated by iteration complexity and computational cost per iteration \cite{bubeck2015convex}.

First, we consider the iteration complexity of (\ref{1}) in the following cases.
\begin{itemize}
  \item If $f_{ij}$ is a simple-structured smooth function (e.g., $f_{ij}$ is quadratic), $\nabla f^{\diamond}_{ij}$ (also $\nabla p_{ij}$ by (\ref{6})) can be obtained efficiently \cite[Sec. 3.3.1]{boyd2004convex}. Widely discussed optimization problems with smooth + nonsmooth cost functions can be referred to in \cite{beck2014fast,hans2009bayesian,zhao2017scope}.
  \item In (\ref{1}), the proximal mapping of $g_{ij}^{\diamond}$ can be efficiently obtained if $g_{ij}$ is simple-structured. For example, consider a regularization problem, where the penalty is a Euclidean $e$-norm: $g_{ij}(\mathbf{x}_{ij})= \| \mathbf{x}_{ij} \|_e$ (e.g., it can be an LASSO problem if $e=1$ \cite{hans2009bayesian}). Then, we can have
\begin{align}
g^{\diamond}_{ij}( \bm{\mu}_{ij}) = \mathds{I}_{Y_{ij}} (\bm{\mu}_{ij}) = \left\{\begin{array}{ll}
                    0 & \hbox{if $\bm{\mu}_{ij} \in Y_{ij}$} \\
                    +\infty  & \hbox{otherwise}
                  \end{array}
                  \right.,
\end{align}
where $Y_{ij}=\{ \mathbf{v}  \in \mathds{R}^{M} | \| \mathbf{v} \|^*_{e} \leq 1 \}$ with $\| \cdot \|^*_{e}$ being the dual norm of $\| \cdot \|_e$. The first equality holds by computing the conjugate of an $e$-norm \cite[Sec. 3.3.1]{boyd2004convex}. Then in (\ref{1}), the proximal mapping of $g^{\diamond}_{ij}$ is a Euclidean projection onto $Y_{ij}$ \cite[Sec. 1.2]{parikh2014proximal}.
\end{itemize}

Based on the above discussion, provided that $f_{ij}$ and $g_{ij}$ are some simple-structured functions, (\ref{1}) only requires some simple operations (e.g., addition, multiplication, and Euclidean projection) with iteration complexity $O(1)$ without any costly inner-loop optimization. Then the overall iteration complexity of (\ref{1}) can be $O(\frac{1}{\epsilon})$ with $\epsilon$ being the convergence error (see Thm. \ref{th1}). Note that for gradient-based iterative algorithms, the computational cost per iteration is linear in dimension \cite{bubeck2015convex}. Then, the overall computational complexity of (\ref{1}) can be $O(\frac{M}{\epsilon})$.

Meanwhile, given that $f_{ij}$ is simple-structured, (\ref{2})-(\ref{5}) can also be computed efficiently only with some simple operations with iteration complexity $O(1)$. Then the overall computational complexity of the CDPG algorithm can be obtained by counting the overall dimension of (\ref{1})-(\ref{5}), which gives $O(\frac{M |\bar{V}|}{\epsilon}) + O(\frac{\sum_{i\in V} n^2_iM}{\epsilon}) + O(\frac{B|\bar{V}|}{\epsilon}) + O(\frac{\sum_{i\in V} \sum_{j\in V_i} |S_{ij}|n_iM}{\epsilon}) + O(\frac{\sum_{i\in V} \sum_{j\in V_i} |\bar{S}_{ij}|B}{\epsilon})$.

\section{Motivating Applications and Numerical Simulation}\label{sa6}

\subsection{Social Welfare Optimization in Commodity Market}

\begin{table*}[htbp]
\caption{Parameters of Simulation A}\label{tm1}
\label{tab2}
\begin{center}
\begin{tabular}{c|cccc|ccc|cc}
\hline
Region $i$ & & \multicolumn{2}{l}{\quad \quad 1} & & &2 & &\multicolumn{2}{l}{\quad \quad3}  \\
\hline
Machine $j$ & 1 & 2 & 3 & 4 & 1 & 2 & 3 & 1 & 2 \\
\hline
$\varpi_{ij}$ & -0.1 & -0.2 & -0.3 & -0.2 & -0.5 & -0.45 & -0.55 & -0.8 & -0.9 \\
\hline
$\varsigma_{ij}$ & 2.1 & 2.2 & 2 & 1.9 & 0.2 & 0.25 & 0.5 & 3.3 & 4.1 \\
\hline
$\underline{x}_{ij}$ & \multicolumn{9}{l}{\quad \quad \quad \quad\quad \quad \quad \quad\quad \quad \quad 0} \\
\hline
$\overline{x}_{ij}$ & 10.5 & 5.5 & 3.33 & 4.75& 0.2 & 0.27 & 0.45 & 2.06 & 2.27 \\
\hline
\end{tabular}
\end{center}
\end{table*}

\begin{table*}[htbp]
\caption{Parameters of Simulation B}\label{tm2}
\label{tab2}
\begin{center}
\begin{tabular}{c|ccccccccc}
\hline
Cluster $i$ & $\underline{x}_i$ & $\overline{x}_i$ & $\alpha_{i,1}$ & $\alpha_{i,2}$ & $\beta_{i,1}$ & $\beta_{i,2}$ & $\rho_{i,1}$ & $\rho_{i,2}$ & $\rho_{i,3}$ \\
\hline
1 & 0.05 & 5 & 100 & 200 & 6.490 & -2.000 & 0.255 & 0.012 & -3.554 \\
\hline
2 & 0.05 & 10 & 120 & 150 & 5.638 & -3.000 & 0.250 & 0.012 & -4.047 \\
\hline
3 & 0.05 & 10 & 40 & 180 & 4.586 & -2.000 & 0.255& 0.012 & -3.094 \\
\hline
\end{tabular}
\end{center}
\end{table*}

In this simulation, we consider a social welfare optimization problem in a commodity market \cite{craven2005optimization}. In this market, we aim to supply certain amount of commodities to multiple consumer regions such that the utility function of the whole consumer community is optimized (the transportation cost is assumed to be negligible).

The utility function of different regions can be obtained by some learning machines based on the regional information \cite{nielsen2004learning}. Due to the possibly large-scale data sets and privacy perseveration issues, it can be inefficient or even infeasible to transmit the whole data sets among the machines. Therefore, distributed learning framework can be employed \cite{verbraeken2020survey}, as introduced as follows.
\begin{itemize}
  \item Distributed learning machines are established, who collect the data in different areas and generate the utility functions based on the local data and learning algorithms.
  \item The utility function of each region can be settled by {\em{ensemble method}} \cite{dietterich2000ensemble}, e.g., take the average of all the generated functions in each region.
\end{itemize}
Then, the machines decide the optimal commodity supply strategy in a distributed manner to optimize the utility function of the whole community.

\begin{figure}[htbp]
  \centering
  \includegraphics[width=7cm]{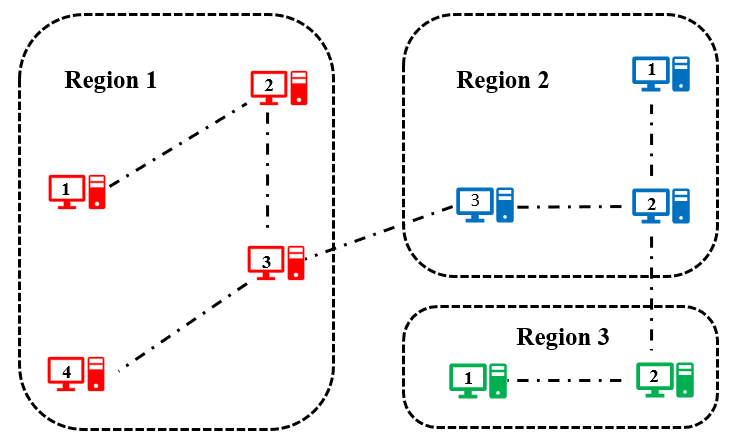}\\
  \caption{Communication typology of learning machines.}\label{mar2}
\end{figure}

Based on the above discussion, we let $f_{ij}(x_i)$ be the utility function generated by the $j$th machine in region $i$. Then, the utility function of region $i$ is settled as $f_i(x_i) = \frac{1}{n_i}\sum_{j\in V_i}f_{ij}(x_i)$ \cite{guo2017learning}. Hence, the social welfare optimization problem can be formulated as
\begin{align}
  \mathrm{(P6)} \quad  \max\limits_{\mathbf{x} \in X} \quad  &  H(\mathbf{x}) \nonumber \\
   \hbox{subject to} \quad & \mathbf{A}\mathbf{x} \leq b. \nonumber
\end{align}
Here, $H(\mathbf{x}) = \sum_{i \in {V}} f_i (x_{i})$, $\mathbf{x}=[x_1,...,x_{|V|}]^{\top}$, $\mathbf{A}= \mathbf{1}_{|V|}^{\top}$, $b$ is total quantity of commodities in store. $X= \prod_{i\in V} \bigcap_{j \in V_i} X_{ij}$ with $X_{ij} = [\underline{x}_{ij},\overline{x}_{ij}]$, where $\underline{x}_{ij}$ and $\overline{x}_{ij}$ are the lower and upper bounds of demand, respectively. In addition, the utility function obtained by the machines is assumed to be in a quadratic form $f_{ij}(x_i) = \varpi_{ij} x^2_i + \varsigma_{ij} x_i$ \cite{craven2005optimization}. $b$ is set as 5. The detailed communication typology of the machines and other parameter settings are shown in Fig. \ref{mar2} and Table \ref{tm1}, respectively.

\begin{figure*}
\centering
\subfigure[Values of $\bm{\lambda}$.]{
\begin{minipage}[t]{0.33\linewidth}
\centering
\includegraphics[width=6.5cm]{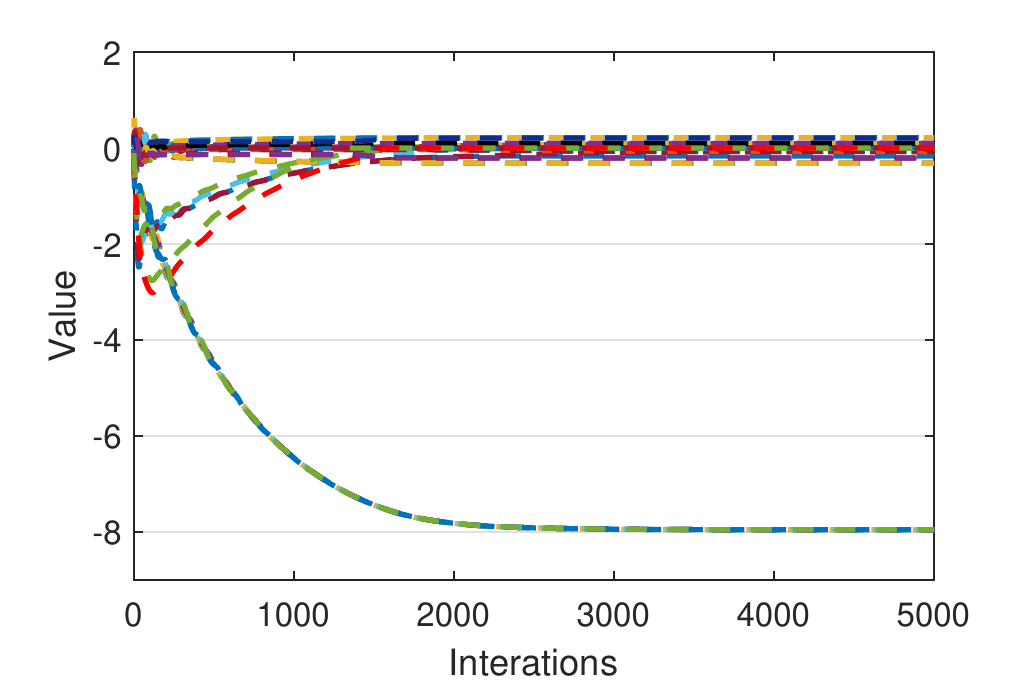}
\end{minipage}%
}%
\subfigure[Values of $\bm{\omega}$.]{
\begin{minipage}[t]{0.33\linewidth}
\centering
\includegraphics[width=6.5cm]{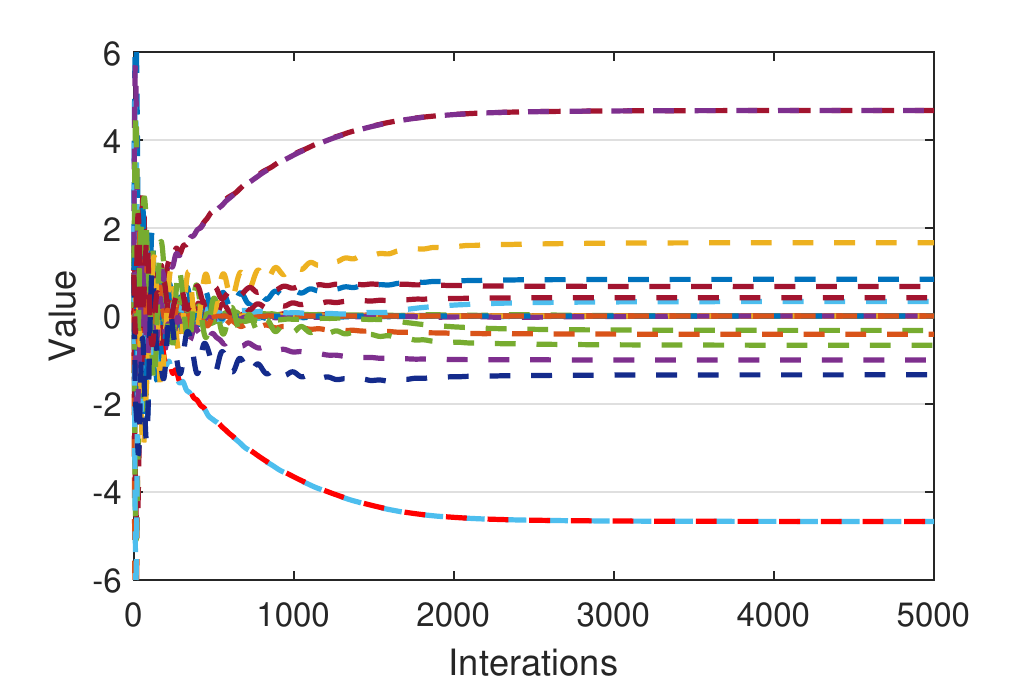}
\end{minipage}%
}%
\subfigure[Values of $o$.]{
\begin{minipage}[t]{0.33\linewidth}
\centering
\includegraphics[width=6.5cm]{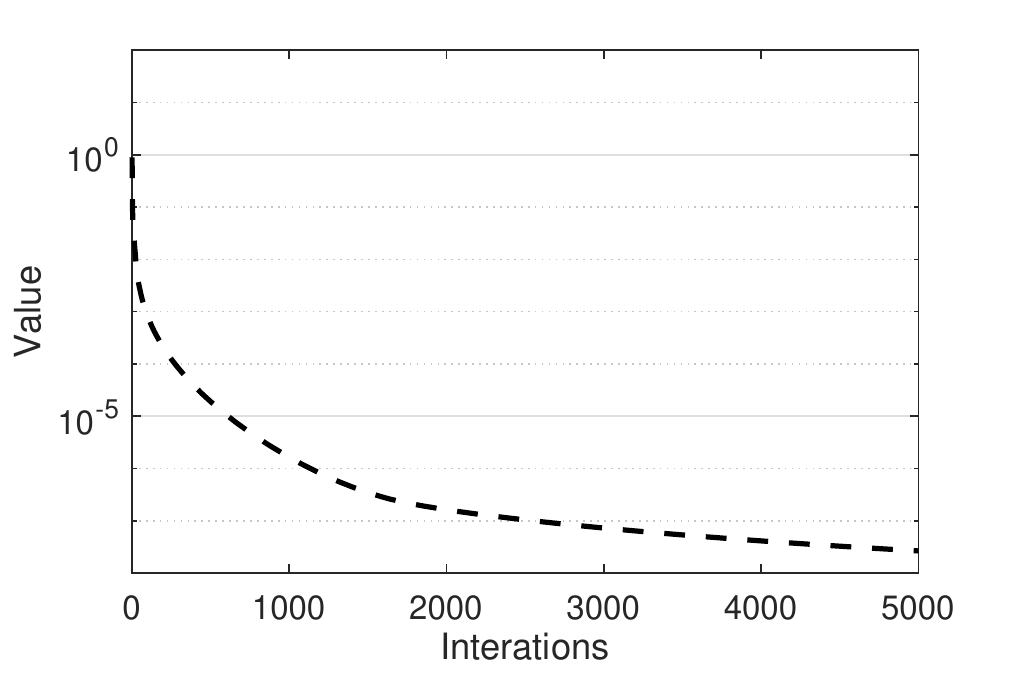}
\end{minipage}%
}%
\caption{Result of Simulation A.}\label{g3}
\end{figure*}

\begin{figure}[htbp]
  \centering
  \includegraphics[width=7cm]{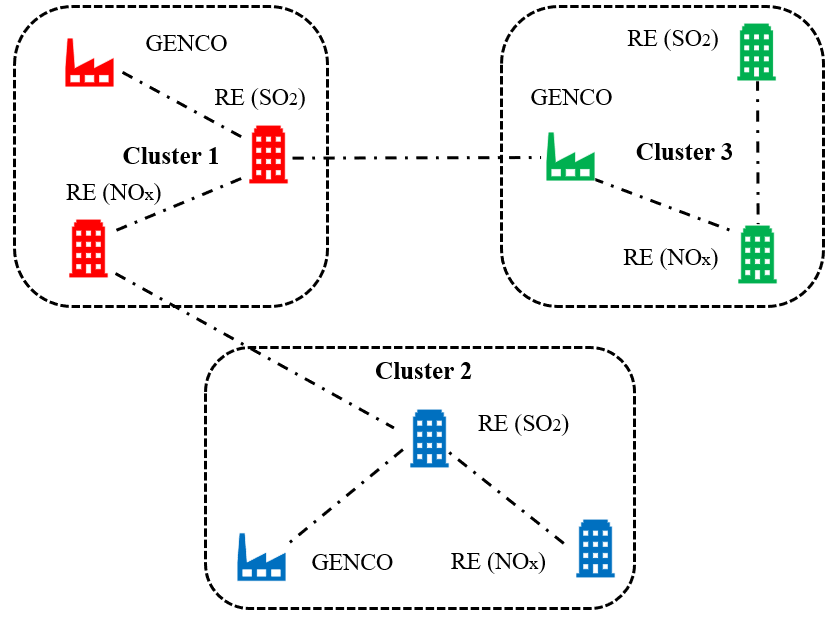}\\
  \caption{Communication typology of GENCOs and REs.}\label{mar}
\end{figure}

\begin{figure*}
\centering
\subfigure[Values of $\bm{\lambda}$.]{
\begin{minipage}[t]{0.33\linewidth}
\centering
\includegraphics[width=6.5cm]{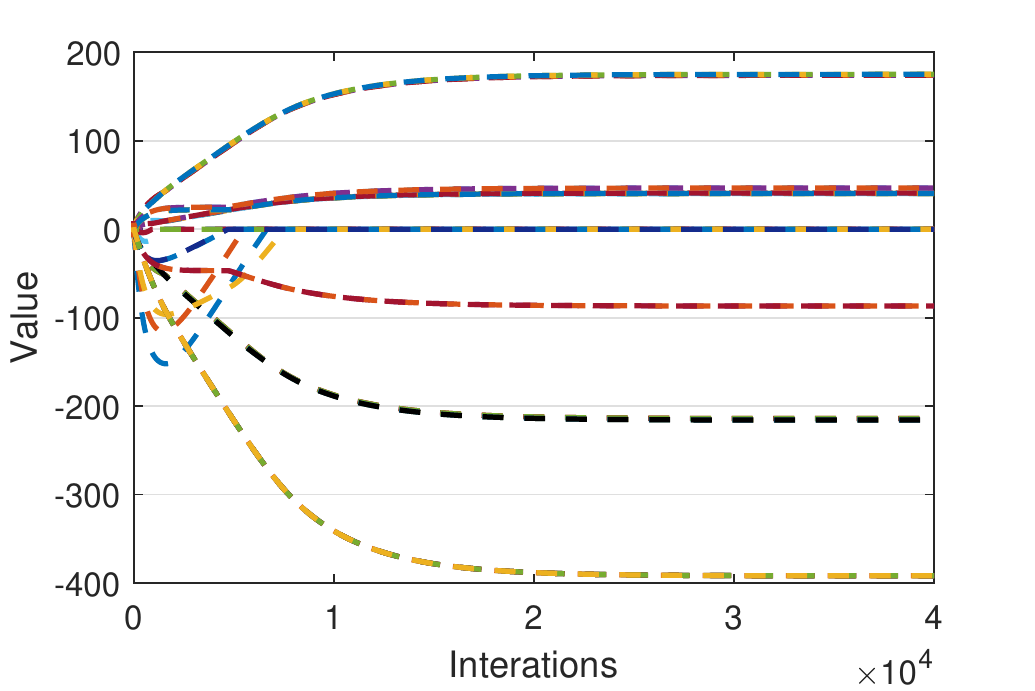}
\end{minipage}%
}%
\subfigure[Values of $\bm{\omega}$.]{
\begin{minipage}[t]{0.33\linewidth}
\centering
\includegraphics[width=6.5cm]{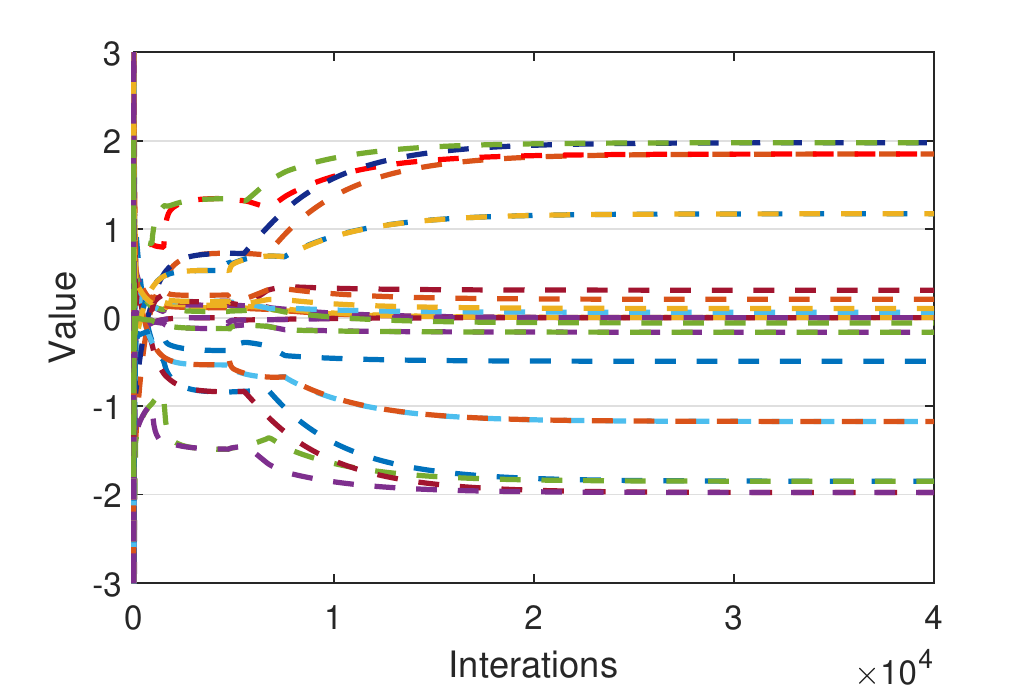}
\end{minipage}%
}%
\subfigure[Values of $o$.]{
\begin{minipage}[t]{0.33\linewidth}
\centering
\includegraphics[width=6.5cm]{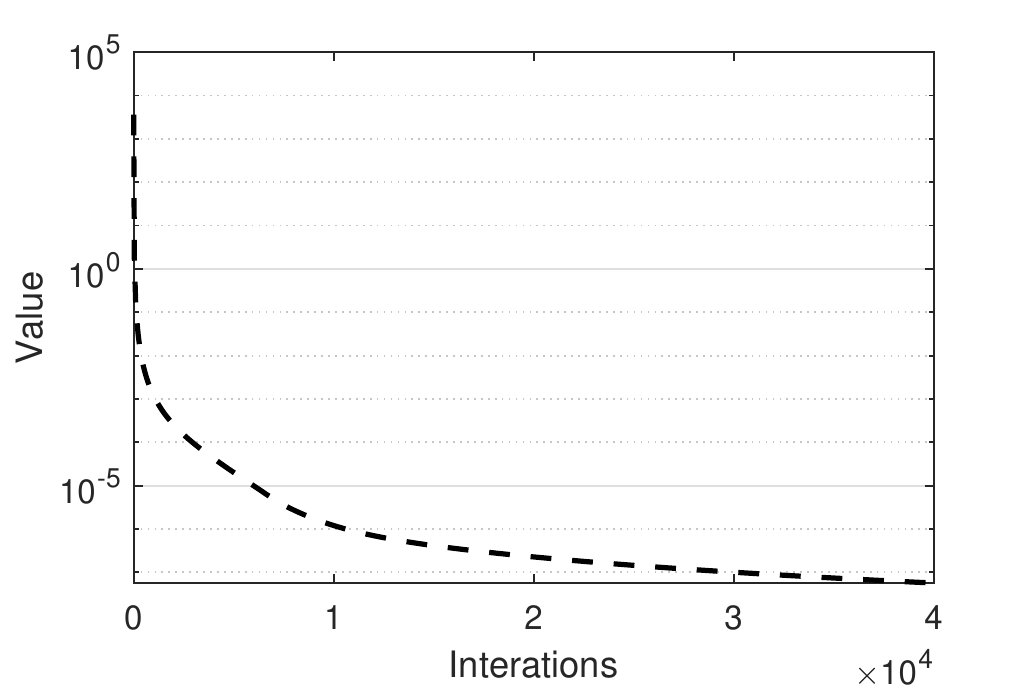}
\end{minipage}%
}%
\caption{Result of Simulation B.}\label{g4}
\end{figure*}

By some direct calculations, the optimal solution to Problem (P6) is $\mathbf{x}^*= [3.33,0,1.67]^{\top}$ (the lower bound of $x_2$ is activated). To characterize the relative convergence error, we define $o=\left|\frac{ \mathcal{L}(\bm{\lambda},\bm{\omega}) - H(\mathbf{x}^*) }{H(\mathbf{x}^*)}\right|$ with certain non-zero $H(\mathbf{x}^*)$. The simulation result is shown in Fig. \ref{g3}. Fig. \ref{g3}-(a) shows that $\bm{\lambda}$ tends to the steady state asymptotically. Fig. \ref{g3}-(b) depicts the trajectory of $\bm{\omega}$. The trajectory of convergence error is shown in Fig. \ref{g3}-(c).

\subsection{Economic Emission Dispatch Optimization in Energy Market}

In this simulation, we consider an economic emission dispatch problem in an energy market. In this market, we consider multiple energy generation companies (GENCOs) and regulation entities (REs), where the latter ones can be some policy makers for the environment's benefit and are responsible for regulating the pollutant emissions \cite{denny2006wind}. Specifically, two sorts of emissions are considered: sulfur dioxide ($\hbox{SO}_2$) and oxides of nitrogen ($\hbox{NO}_{\mathrm{x}}$), which are regulated by possibly different REs. Certain GNECO and its REs are assumed to jointly make decisions to optimize the output of the GNECO. Hence, we treat the GENCO and its REs as a cluster.

In this problem, the objective function of cluster $i$ is designed as $f_i(x_{i})=\chi_iC_i(x_{i})+(1-\chi_i) \delta_i (E^{\mathrm{S}}_i(x_{i})+E_i^{\mathrm{N}}(x_{i}))$ \cite{zou2017new}. Here, $x_i \in X_i = [\underline{x}_i, \overline{x}_i]$ is the quantity of energy generation with $\underline{x}_i$ and $ \overline{x}_i$ being the lower and upper generation limits, respectively. $C_i$, $E^{\mathrm{S}}_i$, and $E^{\mathrm{N}}_i$ are the fuel cost, emission quantity of $\hbox{SO}_2$, and emission quantity of $\hbox{NO}_{\mathrm{x}}$, respectively. $\delta_i>0$ and $\chi_i \in (0,1)$ are the penalty price of emission and weighting factor, respectively. The detailed fuel cost and emission quantity functions are given by
\begin{align}\label{}
& C_i (x_{i})  = \alpha_{i,1} x^2_{i} + \alpha_{i,2} x_{i} + \alpha_{i,3}, \\
& E^{\mathrm{S}}_i (x_{i}) = \beta_{i,1} x^2_{i} + \beta_{i,2} x_{i} + \beta_{i,3}, \\
& E_i^{\mathrm{N}} (x_{i}) = \rho_{i,1} \mathrm{exp} (\rho_{i,2} x_{i}) + \rho_{i,3} x_{i} + \rho_{i,4},
\end{align}
where $\alpha_{i,1}$, $\alpha_{i,2}$, $\alpha_{i,3}$, $\beta_{i,1}$, $\beta_{i,2}$, $\beta_{i,3}$ $\rho_{i,1}$, $\rho_{i,2}$, $ \rho_{i,3}$, and $ \rho_{i,4}$ are parameters. The communication typology of the market is designed in Fig. \ref{mar}. Then, the economic emission dispatch problem of the whole market can be formulated as
\begin{align}
  \mathrm{(P7)} \quad  \min\limits_{\mathbf{x} \in X} \quad  &  \sum_{i \in {V}} f_i (x_{i}) \nonumber \\
   \hbox{subject to} \quad & \mathbf{A}\mathbf{x} = b, \nonumber
   \end{align}
where $\mathbf{A}= \mathbf{1}_{|V|}^{\top}$, $\mathbf{x}=[x_1,...,x_{|V|}]^{\top}$, $X=\prod_{i\in V}X_i$, and $b$ is the total energy demand. The parameters of Problem (P7) are set in Table \ref{tm2} \cite{zou2017new}. Note that $E_i^{\mathrm{N}} (x_{i})$ is strongly convex with a compact $X_i$. The settings of $\alpha_{i,3}$, $\beta_{i,3}$, and $\rho_{i,4}$ are omitted since they are some constants in the objective functions. Without losing the generality, we let $\chi_i=0.5$ and $b=5$. $\delta_i$ is set as $\frac{C_i(\overline{x}_i)}{E^{\mathrm{S}}_i (\overline{x}_{i})+E^{\mathrm{N}}_i (\overline{x}_{i})}$ as suggested by \cite{zou2017new}.

By some direct calculations, the optimal solution to Problem (P7) is $\mathbf{x}^*= [2.38,2.57,0.05]^{\top}$ (the lower bound of $x_3$ is activated). The simulation result is shown in Fig. \ref{g4}. Fig. \ref{g4}-(a) shows that $\bm{\lambda}$ tends to the steady state asymptotically. Fig. \ref{g4}-(b) depicts the trajectory of $\bm{\omega}$. With a similar definition of convergence error as in Simulation A, the trajectory of $o$ is obtained in Fig. \ref{g4}-(c).

\section{Conclusion}\label{sa7}

In this work, we considered a DOP in a multi-cluster network with an affine coupling constraint. In this problem, each cluster can make its decision based on the consensus protocol among the agents involved. To achieve the optimal solution of the whole network, a CDPG algorithm was proposed, where each agent can make updates with local communications. The computational complexity with some simple-structured cost functions were discussed. The performance of the CDPG algorithm was demonstrated by two motivating applications in the simulation.

\appendix{}

\subsection{Proof of Lemma \ref{lam2}}\label{lam2p}

By Lemma \ref{l1}, $ \nabla f^{\diamond}_{ij}$ is Lipschitz continuous with Lipschitz constant $\frac{1}{\sigma_{ij}}$, which means
\begin{align}\label{z1}
& \| \nabla_{\mathbf{v}} f^{\diamond}_{ij}(\mathbf{H}_{ij}\mathbf{v}) - \nabla_{\mathbf{u}} f^{\diamond}_{ij}(\mathbf{H}_{ij}\mathbf{u}) \| \nonumber \\
= & \| \mathbf{H}_{ij}^{\top} \nabla_{\mathbf{H}_{ij}\mathbf{v}} f^{\diamond}_{ij}(\mathbf{H}_{ij}\mathbf{v}) - \mathbf{H}_{ij}^{\top} \nabla_{\mathbf{H}_{ij}\mathbf{u}} f^{\diamond}_{ij}(\mathbf{H}_i\mathbf{u}) \| \nonumber \\
\leq & \| \mathbf{H}_{ij} \| \| \nabla_{\mathbf{H}_{ij}\mathbf{v}} f^{\diamond}_{ij}(\mathbf{H}_{ij}\mathbf{v}) - \nabla_{\mathbf{H}_{ij}\mathbf{u}} f^{\diamond}_{ij}(\mathbf{H}_i\mathbf{u}) \|  \nonumber \\
\leq & \frac{\| \mathbf{H}_{ij} \| }{\sigma_{ij}} \| \mathbf{H}_{ij} \mathbf{v} - \mathbf{H}_{ij} \mathbf{u} \| \nonumber \\
\leq & \frac{\| \mathbf{H}_{ij}  \|^2 }{\sigma_{ij}}\|  \mathbf{v}- \mathbf{u} \|
=  h_{ij} \|  \mathbf{v}- \mathbf{u} \|.
\end{align}
(\ref{z1}) implies $ \nabla_{\bm{\lambda}_{ij}} f^{\diamond}_{ij}(\mathbf{H}_{ij}\bm{\lambda}_{ij})$ is Lipschitz continuous with constant $h_{ij}$, which means $\nabla_{\bm{\lambda}_{ij}} p_{ij} (\bm{\lambda}_{ij})=\nabla_{\bm{\lambda}_{ij}} f^{\diamond}_{ij}(\mathbf{H}_{ij}\bm{\lambda}_{ij}) + \mathbf{E}^{\top}_{ij}$ is also Lipschitz continuous with constant $h_{ij}$.

\subsection{{Proof of Theorem \ref{th1}}}\label{th1p}

By the first-order optimality condition of (\ref{f1}) in terms of (\ref{13}), we have
\begin{align}\label{e-1}
\mathbf{0}  \in & \partial_{\bm{\lambda}} Q (\bm{\lambda}^{t+1}) + \bar{\mathbf{S}}[c](\bm{\lambda}^{t+1} - \bm{\lambda}^t)  + \nabla_{\bm{\lambda}} P (\bm{\lambda}^t) \nonumber \\
&  + \mathbf{Z}^{\top} \mathbf{D}[\pi] \mathbf{Z} \bm{\lambda}^t + \mathbf{Z}^{\top} \bm{\omega}^t  \nonumber \\
 = &  \partial_{\bm{\lambda}} Q (\bm{\lambda}^{t+1}) - \bar{\mathbf{S}}[c](\bm{\lambda}^t- \bm{\lambda}^{t+1}) + \nabla_{\bm{\lambda}} P (\bm{\lambda}^t) \nonumber \\
&  + \mathbf{Z}^{\top} \mathbf{D}[\pi] \mathbf{Z} \bm{\lambda}^t +  \mathbf{Z}^{\top} \bm{\omega}^{t+1}- \mathbf{Z}^{\top} \mathbf{D}[\pi] \mathbf{Z} \bm{\lambda}^{t+1}.
\end{align}
From the convexity of $Q(\bm{\lambda})$, we have
\begin{align}\label{e2}
Q  (\bm{\lambda})  - & Q(\bm{\lambda}^{t+1}) \geq  (\bm{\lambda}-\bm{\lambda}^{t+1})^{\top} \bar{\mathbf{S}}[c] (\bm{\lambda}^t-\bm{\lambda}^{t+1}) \nonumber \\
& - (\bm{\lambda}-\bm{\lambda}^{t+1})^{\top} \nabla_{\bm{\lambda}} P(\bm{\lambda}^t) -(\bm{\lambda}-\bm{\lambda}^{t+1})^{\top} \mathbf{Z}^{\top} \bm{\omega}^{t+1} \nonumber \\
&  + (\bm{\lambda}-\bm{\lambda}^{t+1})^{\top} \mathbf{Z}^{\top} \mathbf{D}[\pi] \mathbf{Z} ( \bm{\lambda}^{t+1} - \bm{\lambda}^t).
\end{align}
From the convexity and Lipschitz continuous differentiability of $p_{ij}$, we have
\begin{align}\label{e3}
 (  \bm{\lambda} - \bm{\lambda}^{t+1})^{\top} & \nabla_{\bm{\lambda}} P (\bm{\lambda}^t)
= \sum_{i \in {V}} \sum_{j \in {V}_i} (\bm{\lambda}_{ij} - \bm{\lambda}_{ij}^t)^{\top} \nabla_{\bm{\lambda}_{ij}} p_{ij} (\bm{\lambda}_{ij}^t) \nonumber \\
& + \sum_{i \in {V}} \sum_{j \in {V}_i}(\bm{\lambda}_{ij}^t - \bm{\lambda}_{ij}^{t+1})^{\top} \nabla_{\bm{\lambda}_{ij}} p_{ij} (\bm{\lambda}_{ij}^t) \nonumber \\
\leq & \sum_{i \in {V}} \sum_{j \in {V}_i} ( p_{ij}(\bm{\lambda}_{ij}) -  p_{ij}(\bm{\lambda}_{ij}^t) )   \nonumber \\
&  +  \sum_{i \in {V}} \sum_{j \in {V}_i} ( p_{ij}(\bm{\lambda}_{ij}^t)-  p_{ij}(\bm{\lambda}_{ij}^{t+1}) ) \nonumber \\
& + \sum_{i \in {V}} \sum_{j \in {V}_i} \frac{h_{ij}}{2} \|  \bm{\lambda}_{ij}^t - \bm{\lambda}_{ij}^{t+1} \|^2 \nonumber  \\
= & P(\bm{\lambda})  -  P(\bm{\lambda}^{t+1})  + \| \bm{\lambda}^t  - \bm{\lambda}^{t+1} \|^2_{\frac{1}{2}\mathbf{S}[h]}.
\end{align}
By (\ref{f2}), we have
\begin{align}\label{abc}
\mathbf{0}  =  \bar{\mathbf{D}}[\pi](\bm{\omega}^t-\bm{\omega}^{t+1}) + \mathbf{Z} \bm{\lambda}^{t+1}.
\end{align}
Therefore, by multiplying the both sides of (\ref{abc}) by $(\bm{\omega}-\bm{\omega}^{t+1})^{\top}$, we have
\begin{align}\label{r1-1}
 (\bm{\omega} & -\bm{\omega}^{t+1})^{\top} \bar{\mathbf{D}}[\pi] (\bm{\omega}^t - \bm{\omega}^{t+1}) \nonumber \\
  & +  (\bm{\omega}-\bm{\omega}^{t+1})^{\top} \mathbf{Z}\bm{\lambda}^{t+1} = 0.
\end{align}
By adding (\ref{e2}) and (\ref{e3}) together from the both sides, we have
\begin{align}\label{r1}
\Phi & (\bm{\lambda}^{t+1}) - \Phi(\bm{\lambda})
 \leq -(\bm{\lambda} -\bm{\lambda}^{t+1})^{\top} \bar{\mathbf{S}}[c] (\bm{\lambda}^t -\bm{\lambda}^{t+1})  \nonumber\\
&+(\bm{\lambda} -\bm{\lambda}^{t+1})^{\top} \mathbf{Z}^{\top} \bm{\omega}^{t+1} +\| \bm{\lambda}^t - \bm{\lambda}^{t+1}\|^2_{\frac{1}{2}\mathbf{S}[h]} \nonumber\\
&+ (\bm{\lambda} -\bm{\lambda}^{t+1})^{\top} \mathbf{Z}^{\top} \mathbf{D}[\pi] \mathbf{Z} (\bm{\lambda}^t- \bm{\lambda}^{t+1})  \nonumber\\
= & - (\bm{\lambda} -\bm{\lambda}^{t+1})^{\top} \bar{\mathbf{S}}[c] (\bm{\lambda}^t -\bm{\lambda}^{t+1}) \nonumber\\
& - (\bm{\omega} -\bm{\omega}^{t+1})^{\top} \bar{\mathbf{D}}[\pi] (\bm{\omega}^t -\bm{\omega}^{t+1})  \nonumber\\
& - (\bm{\omega}-\bm{\omega}^{t+1})^{\top}\mathbf{Z} \bm{\lambda}^{t+1} +(\bm{\omega}^{t+1})^{\top} \mathbf{Z} \bm{\lambda} \nonumber\\
& -(\bm{\omega}^{t+1})^{\top} \mathbf{Z} \bm{\lambda}^{t+1} + \| \bm{\lambda}^t - \bm{\lambda}^{t+1}\|^2_{\frac{1}{2}\mathbf{S}[h] } \nonumber \\
& + (\bm{\lambda} -\bm{\lambda}^{t+1})^{\top} \mathbf{Z}^{\top} \mathbf{D}[\pi] \mathbf{Z} (\bm{\lambda}^t- \bm{\lambda}^{t+1})  \nonumber\\
= & \|  \bm{\lambda} -\bm{\lambda}^t \|^2_{\frac{1}{2}\bar{\mathbf{S}}[c]} - \|  \bm{\lambda} -\bm{\lambda}^{t+1} \|^2_{\frac{1}{2}\bar{\mathbf{S}}[c]} \nonumber\\
&  -  \| \bm{\lambda}^t -\bm{\lambda}^{t+1}\|^2_{\frac{1}{2}\bar{\mathbf{S}}[c]}  + \|  \bm{\omega} -\bm{\omega}^t \|^2_{\frac{1}{2}\bar{\mathbf{D}}[\pi] } \nonumber\\
& - \|  \bm{\omega} -\bm{\omega}^{t+1}\|^2_{\frac{1}{2}\bar{\mathbf{D}}[\pi] } -  \| \bm{\omega}^t -\bm{\omega}^{t+1}\|^2_{\frac{1}{2}\bar{\mathbf{D}}[\pi] }  \nonumber\\
& +(\bm{\omega}^{t+1})^{\top} \mathbf{Z} \bm{\lambda} - \bm{\omega}^{\top} \mathbf{Z} \bm{\lambda}^{t+1} +  \| \bm{\lambda}^t - \bm{\lambda}^{t+1}\|^2_{\frac{1}{2}\mathbf{S}[h]}  \nonumber\\
& - \| \bm{\lambda} -\bm{\lambda}^t\|^2_{\frac{1}{2} \mathbf{Z}^{\top}\mathbf{D}[\pi]\mathbf{Z}} + \| \bm{\lambda} - \bm{\lambda}^{t+1} \|^2_{\frac{1}{2}\mathbf{Z}^{\top}\mathbf{D}[\pi]\mathbf{Z}} \nonumber \\
& +  \|\bm{\lambda}^t -\bm{\lambda}^{t+1}\|^2_{\frac{1}{2}\mathbf{Z}^{\top}\mathbf{D}[\pi]\mathbf{Z}} \nonumber \\
= & \|  \bm{\lambda} -\bm{\lambda}^t \|^2_{\frac{1}{2}\bar{\mathbf{S}}[c]-\frac{1}{2}\mathbf{Z}^{\top}\mathbf{D}[\pi]\mathbf{Z}} \nonumber \\
& - \|  \bm{\lambda} -\bm{\lambda}^{t+1} \|^2_{\frac{1}{2}\bar{\mathbf{S}}[c]-\frac{1}{2}\mathbf{Z}^{\top}\mathbf{D}[\pi]\mathbf{Z}}  \nonumber\\
& -  \| \bm{\lambda}^t -\bm{\lambda}^{t+1} \|^2_{\frac{1}{2}\bar{\mathbf{S}}[c]-\frac{1}{2}\mathbf{S}[h] - \frac{1}{2}\mathbf{Z}^{\top}\mathbf{D}[\pi]\mathbf{Z}} \nonumber \\
& + \|  \bm{\omega} -\bm{\omega}^t\|^2_{\frac{1}{2}\bar{\mathbf{D}}[\pi]} - \|  \bm{\omega} -\bm{\omega}^{t+1} \|^2_{\frac{1}{2}\bar{\mathbf{D}}[\pi]} \nonumber\\
&  - \| \bm{\omega}^t  -\bm{\omega}^{t+1}\|^2_{\frac{1}{2}\bar{\mathbf{D}}[\pi]}  + (\bm{\omega}^{t+1})^{\top} \mathbf{Z} \bm{\lambda} - \bm{\omega}^{\top} \mathbf{Z} \bm{\lambda}^{t+1},
\end{align}
where we use (\ref{r1-1}) in the first equality and the second equality holds with $\mathbf{v}^{\top} \mathbf{u} = \frac{1}{2} (\| \mathbf{v} \|^2 + \| \mathbf{u} \|^2 - \| \mathbf{v} - \mathbf{u} \|^2) $.

Note that (\ref{r1}) holds for all $\bm{\lambda}$ and $\bm{\omega}$. The proof is conducted by discussing the following two scenarios.

{\em{1) Scenario 1:}} If $\mathbf{Z} \bar{\bm{\lambda}}^{T+1} \neq \mathbf{0}$, by letting $\bm{\lambda}=\bm{\lambda}^*$ and $ \bm{\omega} = 2\| \bm{\omega}^*\| \frac{\mathbf{Z} \bar{\bm{\lambda}}^{T+1}}{\| \mathbf{Z} \bar{\bm{\lambda}}^{T+1} \| }$ in (\ref{r1}), we have
\begin{align}\label{r2+1}
 \Phi & (\bm{\lambda}^{t+1}) - \Phi(\bm{\lambda}^*) +  2\| \bm{\omega}^*\| \frac{(\mathbf{Z} \bar{\bm{\lambda}}^{T+1})^{\top}}{\| \mathbf{Z} \bar{\bm{\lambda}}^{T+1} \| } \mathbf{Z} \bm{\lambda}^{t+1}  \nonumber\\
 \leq &  \|  \bm{\lambda}^* -\bm{\lambda}^t \|^2_{\frac{1}{2}\bar{\mathbf{S}}[c]-\frac{1}{2}\mathbf{Z}^{\top}\mathbf{D}[\pi]\mathbf{Z}} \nonumber \\
& - \|  \bm{\lambda}^* -\bm{\lambda}^{t+1} \|^2_{\frac{1}{2}\bar{\mathbf{S}}[c]-\frac{1}{2}\mathbf{Z}^{\top}\mathbf{D}[\pi]\mathbf{Z}}  \nonumber\\
& + \|  2\| \bm{\omega}^*\| \frac{\mathbf{Z} \bar{\bm{\lambda}}^{T+1}}{\| \mathbf{Z} \bar{\bm{\lambda}}^{T+1} \| }  -\bm{\omega}^t \|^2_{\frac{1}{2}\bar{\mathbf{D}}[\pi]} \nonumber\\
& - \|  2\| \bm{\omega}^*\| \frac{\mathbf{Z} \bar{\bm{\lambda}}^{T+1}}{\| \mathbf{Z} \bar{\bm{\lambda}}^{T+1} \| }  -\bm{\omega}^{t+1} \|^2_{\frac{1}{2}\bar{\mathbf{D}}[\pi]},
\end{align}
where $
0< c_{ij} \leq \frac{1}{h_{ij} + \tau_{\mathrm{max}}(\mathbf{Z}^{\top}\mathbf{D}[\pi]\mathbf{Z})}$
is considered such that $\bar{\mathbf{S}}[c]-\mathbf{S}[h] - \mathbf{Z}^{\top}\mathbf{D}[\pi]\mathbf{Z}$ is positive semi-definite in (\ref{r1}). Summing up (\ref{r2+1}) over $t=0,1,...,T$ gives
\begin{align}\label{r4}
  (T & +1)(\Phi(\bar{\bm{\lambda}}^{T+1}) - \Phi(\bm{\lambda}^*) +  2\| \bm{\omega}^* \| \|\mathbf{Z} \bar{\bm{\lambda}}^{T+1} \|)  \nonumber\\
 \leq & \sum_{t=0}^{T} ( \Phi({\bm{\lambda}}^{t+1}) - \Phi(\bm{\lambda}^*) +  2\| \bm{\omega}^* \| \|\mathbf{Z} \bar{\bm{\lambda}}^{T+1} \| )  \nonumber\\
 \leq & \|  2 \| \bm{\omega}^*\| \frac{\mathbf{Z} \bar{\bm{\lambda}}^{T+1}}{\| \mathbf{Z} \bar{\bm{\lambda}}^{T+1} \| } -\bm{\omega}^0 \|^2_{\frac{1}{2}\bar{\mathbf{D}}[\pi]} \nonumber \\
 & + \|  \bm{\lambda}^* -\bm{\lambda}^0 \|^2_{\frac{1}{2}\bar{\mathbf{S}}[c]-\frac{1}{2}\mathbf{Z}^{\top}\mathbf{D}[\pi]\mathbf{Z}} \nonumber\\
\leq & \|  \bm{\omega}^* \|^2_{4\bar{\mathbf{D}}[\pi]} + \| \bm{\omega}^0 \|^2_{\bar{\mathbf{D}}[\pi]} \nonumber\\
& + \|  \bm{\lambda}^* -\bm{\lambda}^0 \|^2_{\frac{1}{2}\bar{\mathbf{S}}[c]-\frac{1}{2}\mathbf{Z}^{\top}\mathbf{D}[\pi]\mathbf{Z}},
\end{align}
where the first inequality is from the convexity of $\Phi$ and the third inequality is from {\em{Cauchy-Schwarz inequality}}. Therefore,
\begin{align}\label{p1+1}
\Phi (\bar{\bm{\lambda}}^{T+1}) - & \Phi(\bm{\lambda}^*)
\leq  \frac{1}{T+1} (\|  \bm{\omega}^* \|^2_{4\bar{\mathbf{D}}[\pi]} + \| \bm{\omega}^0 \|^2_{\bar{\mathbf{D}}[\pi]} \nonumber\\
& + \|  \bm{\lambda}^* -\bm{\lambda}^0 \|^2_{\frac{1}{2}\bar{\mathbf{S}}[c]-\frac{1}{2}\mathbf{Z}^{\top}\mathbf{D}[\pi]\mathbf{Z}}) \nonumber\\
& - 2\| \bm{\omega}^* \| \|\mathbf{Z} \bar{\bm{\lambda}}^{T+1} \|
\leq  \frac{\Theta}{T+1}.
\end{align}
Based on (\ref{fp1}), (\ref{sad}) and (\ref{k2}), we have
\begin{align}\label{pp1}
& \Phi( {\bm{\lambda}} ) - \Phi(\bm{\lambda}^*) +   \bm{\omega}^{*\top}\mathbf{Z}  {\bm{\lambda}}  \geq 0.
\end{align}
Letting $\bm{\lambda}=\bar{\bm{\lambda}}^{T+1}$ in (\ref{pp1}) gives
\begin{align}\label{p1-1}
\Phi(\bar{\bm{\lambda}}^{T+1}) - \Phi(\bm{\lambda}^*)\geq - \| \bm{\omega}^* \| \|\mathbf{Z} \bar{\bm{\lambda}}^{T+1} \|.
\end{align}
By combining the first inequality in (\ref{p1+1}) and (\ref{p1-1}), we have
\begin{align}\label{p2}
& \|  \bm{\omega}^* \| \| \mathbf{Z} \bar{\bm{\lambda}}^{T+1} \| \leq  \frac{\Theta}{T+1}.
\end{align}
By (\ref{p1-1}) and (\ref{p2}), we have
\begin{align}
& \Phi  (\bar{\bm{\lambda}} ^{T+1})  - \Phi(\bm{\lambda}^*)\geq  - \frac{\Theta}{T+1}.  \label{p3}
\end{align}
By combining (\ref{p1+1}), (\ref{p2}) and (\ref{p3}), (\ref{t1}) and (\ref{t2}) are proved.

{\em{2) Scenario 2:}} If $\mathbf{Z} \bar{\bm{\lambda}}^{T+1} = \mathbf{0}$, we let $\bm{\lambda}=\bm{\lambda}^*$ and $\bm{\omega} = 2 \bm{\omega}^*$ in (\ref{r1}), which directly gives
\begin{align}
& \Phi (\bar{\bm{\lambda}} ^{T+1}) - \Phi(\bm{\lambda}^*) \leq \frac{\Theta}{T+1}
\end{align}
by the same derivation process of (\ref{r2+1})-(\ref{p1+1}). Then, considering that $\| \mathbf{Z} \bar{\bm{\lambda}}^{T+1} \| =0$ and $\Phi(\bar{\bm{\lambda}}^{T+1}) - \Phi(\bm{\lambda}^*) \geq - \| \bm{\omega}^* \|  \| \mathbf{Z}  \bar{\bm{\lambda}}^{T+1} \| = 0$, (\ref{t1}) and (\ref{t2}) hold as well.

\bibliographystyle{IEEEtran}

\bibliography{1myref}

\begin{thebibliography}{10}
\providecommand{\url}[1]{#1}
\csname url@samestyle\endcsname
\providecommand{\newblock}{\relax}
\providecommand{\bibinfo}[2]{#2}
\providecommand{\BIBentrySTDinterwordspacing}{\spaceskip=0pt\relax}
\providecommand{\BIBentryALTinterwordstretchfactor}{4}
\providecommand{\BIBentryALTinterwordspacing}{\spaceskip=\fontdimen2\font plus
\BIBentryALTinterwordstretchfactor\fontdimen3\font minus
  \fontdimen4\font\relax}
\providecommand{\BIBforeignlanguage}[2]{{%
\expandafter\ifx\csname l@#1\endcsname\relax
\typeout{** WARNING: IEEEtran.bst: No hyphenation pattern has been}%
\typeout{** loaded for the language `#1'. Using the pattern for}%
\typeout{** the default language instead.}%
\else
\language=\csname l@#1\endcsname
\fi
#2}}
\providecommand{\BIBdecl}{\relax}
\BIBdecl

\bibitem{luo2014provably}
L.~Luo, N.~Chakraborty, and K.~Sycara, ``Provably-good distributed algorithm
  for constrained multi-robot task assignment for grouped tasks,'' \emph{IEEE
  Transactions on Robotics}, vol.~31, no.~1, pp. 19--30, 2014.

\bibitem{lee2017speeding}
K.~Lee, M.~Lam, R.~Pedarsani, D.~Papailiopoulos, and K.~Ramchandran, ``Speeding
  up distributed machine learning using codes,'' \emph{IEEE Transactions on
  Information Theory}, vol.~64, no.~3, pp. 1514--1529, 2017.

\bibitem{bai2017distributed}
L.~Bai, M.~Ye, C.~Sun, and G.~Hu, ``Distributed economic dispatch control via
  saddle point dynamics and consensus algorithms,'' \emph{IEEE Transactions on
  Control Systems Technology}, vol.~27, no.~2, pp. 898--905, 2017.

\bibitem{guo2017distributed}
F.~Guo, C.~Wen, J.~Mao, G.~Li, and Y.-D. Song, ``A distributed hierarchical
  algorithm for multi-cluster constrained optimization,'' \emph{Automatica},
  vol.~77, pp. 230--238, 2017.

\bibitem{shi2019multi}
C.-X. Shi and G.-H. Yang, ``Multi-cluster distributed optimization via random
  sleep strategy,'' \emph{Journal of the Franklin Institute}, vol. 356, no.~10,
  pp. 5353--5377, 2019.

\bibitem{zhou2017multi}
J.~Zhou, C.~Wang, Y.~Li, P.~Wang, C.~Li, P.~Lu, and L.~Mo, ``A multi-objective
  multi-population ant colony optimization for economic emission dispatch
  considering power system security,'' \emph{Applied Mathematical Modelling},
  vol.~45, pp. 684--704, 2017.

\bibitem{notarnicola2016duality}
I.~Notarnicola, M.~Franceschelli, and G.~Notarstefano, ``A duality-based
  approach for distributed min-max optimization with application to demand side
  management,'' in \emph{2016 IEEE 55th Conference on Decision and Control
  (CDC)}.\hskip 1em plus 0.5em minus 0.4em\relax IEEE, 2016, pp. 1877--1882.

\bibitem{chang2016proximal}
T.-H. Chang, ``A proximal dual consensus admm method for multi-agent
  constrained optimization,'' \emph{IEEE Transactions on Signal Processing},
  vol.~64, no.~14, pp. 3719--3734, 2016.

\bibitem{pang2019randomized}
Y.~Pang and G.~Hu, ``Randomized gradient-free distributed optimization methods
  for a multiagent system with unknown cost function,'' \emph{IEEE Transactions
  on Automatic Control}, vol.~65, no.~1, pp. 333--340, 2019.

\bibitem{ning2018distributed}
B.~Ning, Q.-L. Han, and Z.~Zuo, ``Distributed optimization of multiagent
  systems with preserved network connectivity,'' \emph{IEEE Transactions on
  Cybernetics}, vol.~49, no.~11, pp. 3980--3990, 2018.

\bibitem{wang2019distributed}
X.~Wang, J.~Yan, B.~Jin, and W.~Li, ``Distributed and parallel admm for
  structured nonconvex optimization problem,'' \emph{IEEE Transactions on
  Cybernetics}, 2019.

\bibitem{chang2014distributed}
T.-H. Chang, A.~Nedi{\'c}, and A.~Scaglione, ``Distributed constrained
  optimization by consensus-based primal-dual perturbation method,'' \emph{IEEE
  Transactions on Automatic Control}, vol.~59, no.~6, pp. 1524--1538, 2014.

\bibitem{simonetto2016primal}
A.~Simonetto and H.~Jamali-Rad, ``Primal recovery from consensus-based dual
  decomposition for distributed convex optimization,'' \emph{Journal of
  Optimization Theory and Applications}, vol. 168, no.~1, pp. 172--197, 2016.

\bibitem{mosk2010fully}
D.~Mosk-Aoyama, T.~Roughgarden, and D.~Shah, ``Fully distributed algorithms for
  convex optimization problems,'' \emph{SIAM Journal on Optimization}, vol.~20,
  no.~6, pp. 3260--3279, 2010.

\bibitem{necoara2017fully}
I.~Necoara, V.~Nedelcu, D.~Clipici, and L.~Toma, ``On fully distributed dual
  first order methods for convex network optimization,''
  \emph{IFAC-PapersOnLine}, vol.~50, no.~1, pp. 2788--2793, 2017.

\bibitem{wang2021distributed}
J.~Wang and G.~Hu, ``Distributed discrete-time optimization with coupling
  constraints based on dual proximal gradient method in multi-agent networks,''
  \emph{arXiv preprint arXiv:2108.10652}, 2021.

\bibitem{falsone2017dual}
A.~Falsone, K.~Margellos, S.~Garatti, and M.~Prandini, ``Dual decomposition for
  multi-agent distributed optimization with coupling constraints,''
  \emph{Automatica}, vol.~84, pp. 149--158, 2017.

\bibitem{yang2010distributed}
B.~Yang and M.~Johansson, ``Distributed optimization and games: A tutorial
  overview,'' \emph{Networked Control Systems}, pp. 109--148, 2010.

\bibitem{necoara2015linear}
I.~Necoara and V.~Nedelcu, ``On linear convergence of a distributed dual
  gradient algorithm for linearly constrained separable convex problems,''
  \emph{Automatica}, vol.~55, pp. 209--216, 2015.

\bibitem{alghunaim2019proximal}
S.~A. Alghunaim, K.~Yuan, and A.~H. Sayed, ``A proximal diffusion strategy for
  multiagent optimization with sparse affine constraints,'' \emph{IEEE
  Transactions on Automatic Control}, vol.~65, no.~11, pp. 4554--4567, 2019.

\bibitem{li2019gossip}
D.-K. Li, C.-X. Shi, and G.-H. Yang, ``Gossip-based distributed hierarchical
  algorithm for multi-cluster constrained optimisation,'' \emph{IET Control
  Theory \& Applications}, vol.~13, no.~15, pp. 2346--2355, 2019.

\bibitem{li2020distributed}
X.~Li, G.~Feng, and L.~Xie, ``Distributed proximal algorithms for multi-agent
  optimization with coupled inequality constraints,'' \emph{IEEE Transactions
  on Automatic Control}, 2020.

\bibitem{notarnicola2017duality}
I.~Notarnicola and G.~Notarstefano, ``A duality-based approach for distributed
  optimization with coupling constraints,'' \emph{IFAC-PapersOnLine}, vol.~50,
  no.~1, pp. 14\,326--14\,331, 2017.

\bibitem{notarnicola2016asynchronous}
------, ``Asynchronous distributed optimization via randomized dual proximal
  gradient,'' \emph{IEEE Transactions on Automatic Control}, vol.~62, no.~5,
  pp. 2095--2106, 2016.

\bibitem{beck2014fast}
A.~Beck and M.~Teboulle, ``A fast dual proximal gradient algorithm for convex
  minimization and applications,'' \emph{Operations Research Letters}, vol.~42,
  no.~1, pp. 1--6, 2014.

\bibitem{kim2016fast}
D.~Kim and J.~A. Fessler, ``Fast dual proximal gradient algorithms with rate
  $o(1/k^{1.5})$ for convex minimization,'' \emph{arXiv preprint
  arXiv:1609.09441}, 2016.

\bibitem{wang2021composite}
J.~Wang and G.~Hu, ``Composite optimization with coupling constraints via dual
  proximal gradient method with applications to asynchronous networks,''
  \emph{arXiv preprint arXiv:2102.12797}, 2021.

\bibitem{borwein2010convex}
J.~Borwein and A.~S. Lewis, \emph{Convex analysis and nonlinear optimization:
  theory and examples}.\hskip 1em plus 0.5em minus 0.4em\relax Springer Science
  \& Business Media, 2010.

\bibitem{beck2017first}
A.~Beck, \emph{First-order methods in optimization}.\hskip 1em plus 0.5em minus
  0.4em\relax SIAM, 2017.

\bibitem{chung1997spectral}
F.~R. Chung and F.~C. Graham, \emph{Spectral graph theory}.\hskip 1em plus
  0.5em minus 0.4em\relax American Mathematical Soc., 1997, no.~92.

\bibitem{dimarogonas2010stability}
D.~V. Dimarogonas and K.~H. Johansson, ``Stability analysis for multi-agent
  systems using the incidence matrix: Quantized communication and formation
  control,'' \emph{Automatica}, vol.~46, no.~4, pp. 695--700, 2010.

\bibitem{hans2009bayesian}
C.~Hans, ``Bayesian lasso regression,'' \emph{Biometrika}, vol.~96, no.~4, pp.
  835--845, 2009.

\bibitem{zhao2017scope}
S.-Y. Zhao, R.~Xiang, Y.-H. Shi, P.~Gao, and W.-J. Li, ``Scope: Scalable
  composite optimization for learning on spark,'' in \emph{Thirty-First AAAI
  Conference on Artificial Intelligence}, 2017.

\bibitem{shi2015proximal}
W.~Shi, Q.~Ling, G.~Wu, and W.~Yin, ``A proximal gradient algorithm for
  decentralized composite optimization,'' \emph{IEEE Transactions on Signal
  Processing}, vol.~63, no.~22, pp. 6013--6023, 2015.

\bibitem{schmidt2011convergence}
M.~Schmidt, N.~Roux, and F.~Bach, ``Convergence rates of inexact
  proximal-gradient methods for convex optimization,'' \emph{Advances in Neural
  Information Processing Systems}, vol.~24, pp. 1458--1466, 2011.

\bibitem{chang2014multi}
T.-H. Chang, M.~Hong, and X.~Wang, ``Multi-agent distributed optimization via
  inexact consensus admm,'' \emph{IEEE Transactions on Signal Processing},
  vol.~63, no.~2, pp. 482--497, 2014.

\bibitem{florea2020generalized}
M.~I. Florea and S.~A. Vorobyov, ``A generalized accelerated composite gradient
  method: Uniting nesterov's fast gradient method and fista,'' \emph{IEEE
  Transactions on Signal Processing}, vol.~68, pp. 3033--3048, 2020.

\bibitem{boyd2004convex}
S.~Boyd, S.~P. Boyd, and L.~Vandenberghe, \emph{Convex optimization}.\hskip 1em
  plus 0.5em minus 0.4em\relax Cambridge University Press, 2004.

\bibitem{rockafellar1970convex}
R.~T. Rockafellar, \emph{Convex analysis}.\hskip 1em plus 0.5em minus
  0.4em\relax Princeton university press, 1970, no.~28.

\bibitem{hanson1981sufficiency}
M.~A. Hanson, ``On sufficiency of the kuhn-tucker conditions,'' \emph{Journal
  of Mathematical Analysis and Applications}, vol.~80, no.~2, pp. 545--550,
  1981.

\bibitem{parikh2014proximal}
N.~Parikh and S.~Boyd, ``Proximal algorithms,'' \emph{Foundations and Trends in
  Optimization}, vol.~1, no.~3, pp. 127--239, 2014.

\bibitem{bubeck2015convex}
S.~Bubeck, ``Convex optimization: Algorithms and complexity,''
  \emph{Foundations and Trends{\textregistered} in Machine Learning}, vol.~8,
  no. 3-4, pp. 231--357, 2015.

\bibitem{craven2005optimization}
B.~D. Craven and S.~M. Islam, \emph{Optimization in economics and finance: some
  advances in non-linear, dynamic, multi-criteria and stochastic models}.\hskip
  1em plus 0.5em minus 0.4em\relax Springer Science \& Business Media, 2005,
  vol.~7.

\bibitem{nielsen2004learning}
T.~D. Nielsen and F.~V. Jensen, ``Learning a decision maker's utility function
  from (possibly) inconsistent behavior,'' \emph{Artificial Intelligence}, vol.
  160, no. 1-2, pp. 53--78, 2004.

\bibitem{verbraeken2020survey}
J.~Verbraeken, M.~Wolting, J.~Katzy, J.~Kloppenburg, T.~Verbelen, and J.~S.
  Rellermeyer, ``A survey on distributed machine learning,'' \emph{ACM
  Computing Surveys (CSUR)}, vol.~53, no.~2, p. 3377454, 2020.

\bibitem{dietterich2000ensemble}
T.~G. Dietterich, ``Ensemble methods in machine learning,'' in
  \emph{International workshop on multiple classifier systems}.\hskip 1em plus
  0.5em minus 0.4em\relax Springer, 2000, pp. 1--15.

\bibitem{guo2017learning}
Z.-C. Guo, L.~Shi, and Q.~Wu, ``Learning theory of distributed regression with
  bias corrected regularization kernel network,'' \emph{The Journal of Machine
  Learning Research}, vol.~18, no.~1, pp. 4237--4261, 2017.

\bibitem{denny2006wind}
E.~Denny and M.~O'Malley, ``Wind generation, power system operation, and
  emissions reduction,'' \emph{IEEE Transactions on Power Systems}, vol.~21,
  no.~1, pp. 341--347, 2006.

\bibitem{zou2017new}
D.~Zou, S.~Li, Z.~Li, and X.~Kong, ``A new global particle swarm optimization
  for the economic emission dispatch with or without transmission losses,''
  \emph{Energy Conversion and Management}, vol. 139, pp. 45--70, 2017.

\end{thebibliography}

\end{document}